\documentclass{article}
\usepackage[utf8]{inputenc}
\usepackage{graphicx} 
\usepackage{amsmath}
\usepackage{amssymb}
\usepackage{booktabs}
\usepackage{siunitx}
\usepackage{algorithm}
\usepackage{fullpage}
\usepackage{algpseudocode}
\usepackage[numbers]{natbib}
\usepackage{tikz}
\usetikzlibrary{trees}
\usepackage[T1]{fontenc}
\usepackage{authblk}
\usepackage[colorlinks=true, citecolor=blue, linkcolor=blue, urlcolor=cyan]{hyperref}

\usepackage{tikz}
\usetikzlibrary{shapes.geometric, arrows, positioning, calc}

\tikzset{
  startstop/.style={
    rectangle,
    rounded corners,
    minimum width=3cm,
    minimum height=1cm,
    text centered,
    draw=black,
    fill=blue!20
  },
  process/.style={
    rectangle,
    minimum width=3cm,
    minimum height=1cm,
    text centered,
    draw=black,
    fill=orange!20
  },
  processBoth/.style={
    process,
    fill=none
  },
  process0/.style={
    process,
    fill=red!20
  },
  process1/.style={
    process,
    fill=yellow!20
  },
  decision/.style={
    diamond,
    aspect=2,
    draw=black,
    fill=purple!20,
    text centered,
    inner sep=1pt
  },
  arrow/.style={
    thick,
    ->,
    >=stealth
  }
}
\tikzstyle{arrow} = [thick,->,>=stealth]

\usepackage{ulem}
\usepackage{xcolor}
\definecolor{myrvs}{rgb}{1.0,0.01,0.24}

\title{A multi-mesh adaptive finite element method for solving the Gross-Pitaevskii equation}
\author[1]{Mingzhe Li}
\author[1]{Yang Kuang\thanks{Corresponding author. Email: \texttt{ykuang@gdut.edu.cn}}}
\author[2,3]{Zhicheng Hu}
\affil[1]{School of Mathematics and Statistics \& Center for Mathematics and Interdisciplinary Science (CMIS), Guangdong University of Technology,  China}
\affil[2]{School of Mathematics, Nanjing University of Aeronautics and Astronautics, Nanjing, 211106, China}
\affil[3]{Key Laboratory of Mathematical Modelling and High Performance Computing of Air Vehicles (NUAA), MIIT, Nanjing, 211106, China}
\date{\today}

\usepackage{etoolbox}
\makeatletter
\patchcmd{\@afterheading}{\@nobreaktrue}{\@nobreakfalse}{}{}
\makeatother
\newenvironment{keywords}{%
    \par\bigskip\noindent\textbf{{Keywords:}}\quad\ignorespaces
}{
}

\begin{document}

    \maketitle
    \begin{abstract}
        It is found that the wave functions of the Gross-Pitaevskii equation (GPE) often vary significantly in different spatial regions, with some components exhibiting sharp variations while others remain smooth. Solving the GPE on a single mesh, even with adaptive refinement, can lead to excessive computational costs due to the need to accommodate the most oscillatory solution. To address this issue, we present a multi-mesh adaptive finite element method for solving the GPE. To this end, we first convert it into a time-dependent equation through the imaginary time propagation method. Then the equation is discretized by the backward Euler method temporally and the multi-mesh adaptive finite element method spatially. The proposed method is compared with the single-mesh adaptive method through a series of numerical experiments, which demonstrate that the multi-mesh adaptive method can achieve the same numerical accuracy with less computational consumption.
    \end{abstract}
            

    \section{Introduction}
    A Bose–Einstein condensate (BEC) represents an extreme quantum state of matter formed when a dilute bosonic gas is cooled to temperatures approaching absolute zero. The condensates of boson particles were first predicted by S.N. Bose and A. Einstein in 1924 \cite{Bose1924,Einstein1924}, and experimentally achieved under ultra-low temperatures in 1995 \cite{Anderson1995,Davis1995}. When the temperature $\mathit{T}$ is much less than the critical condensation temperature $\mathit{T_c}$, the BEC can be accurately described by a wave function $\psi(\mathbf{x},\mathit{t})$, which is governed by the nonlinear Schr\"odinger equation, i.e., the Gross-Pitaevskii equation (GPE) \cite{Gross1961,Pitaevskii1961}.

    Over the past two decades, a wide range of numerical methods have been developed for the efficient and accurate solution of the GPE. Among early studies, the introduction of sinusoidal pseudospectral schemes based on the normalized gradient flow, specifically the backward Euler sine pseudospectral and backward/forward Euler sine-pseudospectral method \cite{Bao2006,Bao2013}, provided spatial accuracy and temporal stability, while guaranteeing monotonic energy decay in computing both ground and first excited states of Bose-Einstein condensates (BECs). Several alternative gradient flows have been developed by modifying the underlying geometry of the gradient flow, including the projected Sobolev gradient flow \cite{Henning2020,Chen2024} and the J-method \cite{Altmann2021}, both of which have been successfully applied to compute ground and excited states of the BEC. Apart from the gradient flow based methods, the ground state BEC was reformulated as an energy minimization problem. The demand of higher accuracy in dynamical simulations motivated the development of the time-splitting spectral method \cite{Antoine2013}, which preserves key physical invariants and exhibits excellent scalability in high-dimensional computations. Building upon these ideas, subsequent studies including compact finite difference method \cite{Zhang2015}, discontinuous Galerkin finite element method \cite{Yang2025}, Krylov solvers \cite{Zhang2016}, Newton-based method \cite{Huang2022}, achieved further gains in both accuracy and efficiency for solving the GPE.

    Given its potential applications, a core aim of BEC research is to form condensates with as many particles as possible \cite{Bao2013}. It is particularly noteworthy to study a scenario where two (or more) independently growing condensate are merged into a large single condensate to achieve this goal. Consequently, efficient and accurate solving multi-component Gross-Pitaevskii equations (GPEs) is important.

    For solving multi-component GPEs, the coupling between condensate components introduces additional computational challenges, prompting the development of a range of numerical methods. Early studies combined the Thomas-Fermi approximation with geometric coordinate transformations to reduce the two-component GPEs to a piecewise spherical-parabola assembly problem, enabling the construction of general topological phase \cite{Riboli2002}. Building on the framework of normalized gradient flow, a unified and stable computational scheme based on backward Euler and time-splitting sine-spectral method was established, ensuring energy and particle number conservation across arbitrary dimensions, component numbers, and interaction strengths \cite{Bao2004}, and extensions to systems with internal Josephson coupling established the existence, uniqueness, and asymptotic properties of ground states through variational theory and normalized gradient-flow algorithm \cite{Bao2011}. Subsequently, through a combination of the Thomas-Fermi approximation and variational analysis, investigated the spatial localization and phase segregation of ground and excited states in two-dimensional two-component GPEs under strong interspecies repulsion \cite{Caliari2008} . Later, a spectral method based on energy functional minimization was proposed, which expands the wave functions into Hermite basis functions and applies a simplified Newton optimization algorithm \cite{Caliari2009}. This method effectively computes ground states of multi-component BECs, particularly improving computational efficiency and stability in high-dimensional and strongly coupled regimes, while accurately capturing phase-separation phenomena. Recently, a time-dependent external harmonic potential was employed to investigate non-autonomous vector solitons in multi-component GPEs, offering a description of the two-component GPEs corresponding to two hyperfine states \cite{Liu2018}. Nevertheless, the accurate and efficient computation of ground and excited states under complex multi-component interactions remains an open challenge, underscoring the need for efficient and accurate numerical methods for multi-component GPEs.

    Condensate wave functions governed by the GPE often exhibit strong spatial variations or even singular behavior, particularly near trapping centers or in the presence of nonlinear interactions. Such localized features significantly reduce the solution regularity, implying that large uniform meshes are required to achieve accuracy, thus leading to prohibitive computational cost. $h$-Adaptive finite element methods \cite{Li2005,Chen2025,Tao2022,Verfurth2018} address this challenge by providing locally refined approximation spaces that efficiently capture high-gradient regions while coarsening elsewhere. In many GPE applications, one must compute multiple solutions—including ground and excited states, or multi-component condensates—whose regularity and singular structures can differ substantially. If all components are solved on a single adaptive mesh, the refinement must accommodate the most irregular solution, causing unnecessary degrees of freedom for the others. This motivates constructing multiple approximation spaces that adapt independently to the solution behavior of each component. The multi-mesh adaptive finite element method proposed in this work is designed precisely to achieve this improved balance between accuracy and computational complexity.

    The use of multiple meshes, while offering considerable potential for efficiency, introduces additional algorithmic and numerical challenges beyond those encountered in single-mesh approaches. In particular, (i) the simultaneous management of several independently refined meshes requires a robust strategy for local $h$-adaptation, and (ii) accurate multi-mesh communication must be ensured, especially in the evaluation of nonlinear and coupling terms involving integrals over the physical domain. To resolve these issues, we adopt the hierarchical geometry tree data structure \cite{Li2005}, which naturally supports localized refinement and coarsening while maintaining efficient mesh organization. Additionally, the inter-mesh communication is realized through a quadrature-point-based interpolation strategy that avoids projection, thereby preventing degradation of accuracy in coupled computations. As a result, the proposed methodology enables each component of the GPEs to be approximated in an optimal finite element space tailored to its own physical characteristics, achieving high-fidelity simulations with significantly reduced computational expense compared to conventional uniform-mesh or single-mesh adaptive schemes.

    In virtue of the multi-mesh approach, we establish a multi-mesh adaptive finite element framework to solve the GPE. By solving two wave functions on two different approximation spaces separately, the systematical convergence of total energy and total chemical potential while reducing computational cost relative to the single-mesh method at comparable accuracy. A jump-type \textit{a posteriori}a error estimator is employed to guide mesh adaptation. The resulting numerical accuracy is illustrated by representative numerical examples. By comparing the computational costs of the single-mesh and multi-mesh formulations, the efficiency of the proposed multi-mesh adaptive scheme is demonstrated.

    In this work, we organize the paper as follows. We first introduce the GPE, numerical discretization, and adaptive methods in Section 2. In Section 3,  the multi-mesh adaptive finite element method is demonstrated and its associated numerical challenges are discussed. Section 4 proposes various numerical experiments to verify the effectiveness and robustness of the proposed method. Finally, we conclude this manuscript with a summary of our work.

    \section{The Gross-Pitaevskii equation and finite element discretization}
    This section begins with an introduction to the GPE and a discussion of appropriate initial guesses for computing the ground state. Then the discretizations for GPE including the finite element discretization in space and the backward Euler scheme in time are presented. In the end of this section, we introduce the adaptive finite element method to accelerate the computation.

    \subsection{Gross-Pitaevskii equation}
    \label{subsec:2.1}
    The $d$-dimensional GPE is given by
    \begin{equation}
        \left\{
        \begin{aligned}
        \label{eq:2}
        & \mathrm{i}\frac{\partial}{\partial t}\psi(\mathbf{x},t)  =-\frac{1}{2}\Delta\psi(\mathbf{x},t)+V(\mathbf{x})\psi(\mathbf{x},t)+\beta|\psi(\mathbf{x},t)|^2\psi(\mathbf{x},t),\quad \mathbf{x}\in\Omega\subseteq\mathbb{R}^d\left(d=1,2,3,\ldots\right),\quad t>0, 
        \\
        & \psi(\mathbf{x},t) = 0,\quad \mathbf{x}\in\partial\Omega,\quad t\geqslant0,
    \end{aligned}\right.
    \end{equation}
    where $\mathrm{i}$ is the imaginary unit, $t$ is time, the parameter $\beta$ characterizes the interaction among atoms, $V(x)$ is the trapping potential and is usually harmonic taking the following form:
    \begin{equation}
        \label{eq:3}
        V(\mathbf{x}) = \frac{1}{2} \left( \gamma_1^2 x_1^2 + \ldots + \gamma_d^2 x_d^2 \right),
    \end{equation}
    where $\mathbf{x} = {(x_{1},x_{2},...,x_{d})}^T\in\mathbb{R}^d$ and $\gamma_{j} > 0$ for $1 \leq j \leq d$.
    The normalization of the wave function is
    \begin{equation}
        \label{eq:4}
        N(\psi)=\|\psi(\cdot,t)\|^2=\int_{\Omega}|\psi(\mathbf{x},t)|^2d\mathbf{x}=1,\quad t>0,
    \end{equation}
    and the energy per particle is written as
    \begin{equation}
        \label{eq:5}
        E(\psi(\cdot,t))=\int_{\Omega}\left[\frac{1}{2}|\nabla\psi(\cdot,t)|^2+V(\mathbf{x})|\psi(\cdot,t)|^2+\frac{\beta}{2}|\psi(\cdot,t)|^4\right]d\mathbf{x},\quad t>0.
    \end{equation}
    The formal solution of Eq. \eqref{eq:2} takes the following form:
    $$\psi(\mathbf{x},t)=e^{-\mathrm{i}\mu t}\psi(\mathbf{x}),$$
    where 
    $\psi(\mathbf{x})$ denotes the time-independent wave function and $\mu$ stands for the chemical potential. Substituting this expression into Eq. \eqref{eq:2} leads to the following eigenvalue problem for $\mu$:

    \begin{equation}\left\{\begin{aligned}
        \label{eq:6}
        \mu\psi(\mathbf{x}) & =-\frac{1}{2}\Delta\psi(\mathbf{x})+V(\mathbf{x})\psi(\mathbf{x})+\beta|\psi(\mathbf{x})|^{2}\psi(\mathbf{x}),\quad \mathbf{x}\in\Omega, 
        \\
        \psi(\mathbf{x}) & =0,\quad \mathbf{x}\in\partial\Omega.
    \end{aligned}\right.
    \end{equation}
    The normalization of $\psi$ is as follows:
    \begin{equation}
        \label{eq:7}
        \left\|\psi\right\|^2=\int_\Omega\left|\psi(\mathbf{x})\right|^2d\mathbf{x}=1.
    \end{equation}
    Eq. \eqref{eq:6} defines a nonlinear eigenvalue problem, in which the eigenvalue  $\mu$ can be evaluated from the associated eigenfunction $\psi(\mathbf{x})$ as
    \begin{equation}\begin{aligned}
        \label{eq:8}
        \mu & =\int_{\Omega}\left[\frac{1}{2}|\nabla\psi(\mathbf{x})|^{2}+V(\mathbf{x})|\psi(\mathbf{x})|^{2}+\beta|\psi(\mathbf{x})|^{4}\right]d\mathbf{x} \\
        & =E(\psi)+\int_{\Omega}\frac{\beta}{2}|\psi(\mathbf{x})|^{4}d\mathbf{x}.
    \end{aligned}\end{equation}

    In this work, instead of solving the eigenvalue problem directly, we apply the imaginary time propagation method \cite{Chiofalo2000,Chin2009,Lehtovaara2007,Satarić2016,YoungS2010,Kuang2018} where an imaginary time variable $t:=-\mathrm{i}t$ is introduced. Thus Eq. \eqref{eq:2} is transformed into the following time-dependent equation
    \begin{equation}\begin{aligned}
        \label{eq:9}
        \frac{\partial}{\partial t}\psi(\mathbf{x},t) & =\frac{1}{2}\Delta\psi(\mathbf{x},t)-V(\mathbf{x})\psi(\mathbf{x},t)-\beta\left|\psi\right|^2\psi(\mathbf{x},t),\quad \mathbf{x}\in\Omega,\quad t_0\leqslant t<t_{n+1}, 
        \\
        \psi(\mathbf{x},t) & =0,\quad \mathbf{x}\in\partial\Omega, 
        \\
        \psi(\mathbf{x},0) & =\psi^0(\mathbf{x}),\quad \mathbf{x}\in\Omega.
    \end{aligned}\end{equation}
    It can be derived that the solution of Eq. \eqref{eq:9} converges to the ground state solution of Eq. \eqref{eq:6} as $t\rightarrow+\infty$ \cite{Lehtovaara2007, Kuang2018}. Eq. \eqref{eq:9} represents a dissipative nonlinear parabolic partial differential equation, for which mass conservation is not preserved by the solution \cite{Yang2025}. Consequently, the solution is projected onto the unit sphere at the end of each time step to enforce this property:
    \begin{equation}
        \label{eq:10}
        {\psi}(\mathbf{x}, t_{n+1}) \, := \, \frac{\psi(\mathbf{x}, t_{n+1})}{\|\psi(\mathbf{x}, t_{n+1})\|}, \,\quad \mathbf{x} \in \Omega.
    \end{equation}

    It is known that a proper initial guess is one of the key issues for efficiently computing the ground state \cite{Bao2006}. Without lose of generality, we assume  the trapping potential $V(x)$ satisfying 
    \begin{equation}
        \label{eq:11}
        V(\mathbf{x})=V_{0}(\mathbf{x})+W(\mathbf{x}),\quad\quad V_{0}(\mathbf{x})=\frac{1}{2}(\gamma_{1}^{2}x_{1}^{2}+...+\gamma_{d}^{2}x_{d}^{2}),
    \end{equation}
    where $\mathbf{x} = (x_{1},x_{2},...,x_{d})^T\in \mathbb{R}^d$ and $\gamma_{j} > 0$ for $1 \leq j \leq d$. A typical example for $W(\mathbf{x})$ is the optical lattice potential \cite{Pitaevskii2003,Bao2007} which is written as
    \begin{equation}
        \label{eq:12}
        W(\mathbf{x})=k_{1}\sin^{2}(q_{1}x_{1}^{2})+...+k_{d}\sin^{2}(q_{d}x_{d}^{2}),
    \end{equation}
    with $k_j$ and $q_j$ being positive constants. When $|\beta|$ is not small, a possible choice of the initial data could be the Thomas-Fermi approximation \cite{Bao2003,Bao2007} as follows:
    \begin{equation}
        \label{eq:13}
        \psi(\mathbf{x,0}) = \frac{\psi_{\mathrm{g}}^{\mathrm{TF}}(\mathbf{x})}{\|\psi_{\mathrm{g}}^{\mathrm{TF}}\|}, \quad
        \psi_{\mathrm{g}}^{\mathrm{TF}}(\mathbf{x}) = 
        \begin{cases} 
            \sqrt{\frac{\mu_{\mathrm{g}}^{\mathrm{TF}} - V(\mathbf{x})}{\beta}}, & V_0(\mathbf{x}) < \mu_{\mathrm{g}}^{\mathrm{TF}}, \\
            0, & \text{otherwise},
        \end{cases}
        \quad \mathbf{x} \in \mathbb{R}^d,
    \end{equation}
    where
    \begin{equation}
        \label{eq:14}
        \mu_{\mathrm{g}}^{\mathrm{TF}} = \frac{1}{2}
        \begin{cases}
            \left(3 \beta \gamma_1 \right)^{2/3}, & d = 1, \\
            \left(4 \beta \gamma_1 \gamma_2 \right)^{1/2}, & d = 2, \\
            \left(15 \beta \gamma_1 \gamma_2 \gamma_3 / (4 \pi) \right)^{2/3}, & d = 3.
        \end{cases}
    \end{equation}

    \subsection{Numerical discretization}
    Owing to the numerical stability of the backward Euler method, we apply this method to discretize the Eq. \eqref{eq:9} on temporal domain. Let the superscript $n$ denote the corresponding term at time $t=n\Delta t$, where $\Delta t$ is the time step. By applying the backward Euler method in Eq. \eqref{eq:9}, we obtain
    \begin{equation}
        \label{eq:15}
        -\frac{1}{\Delta t} \psi^{(n)} = \frac{1}{2} \Delta \psi^{(n+1)} - \left( V(\mathbf{x}) + \beta |\psi^{(n)}|^2 + \frac{1}{\Delta t} \right) \psi^{(n+1)}.
    \end{equation}
    Then we consider the finite element method for the spatial discretization. Specifically, the variational problem can be formulated as: find $\psi\in{V}$ such that 
    \begin{equation}
        \label{eq:16}
        \int_{\Omega} \left[ \frac{1}{2} \nabla \psi^{(n+1)} \nabla v + \left( V(\mathbf{x}) + \beta |\psi^{(n)}|^2 + \frac{1}{\Delta t} \right) \psi^{(n+1)} v \right] d\mathbf{x} = \int_{\Omega} \frac{1}{\Delta t} \psi^{(n)} v \, d\mathbf{x}, \quad \forall v \in  V,
    \end{equation}
    where $V  = \{{v\in{H^1}}:v|_{\mathcal{K}}\in{P^2},\forall \mathcal{K}\in\mathcal{T}_h\}$. The discrete problem of the variational problem mentioned above is: find $\psi\in{V_h}$ such that 
    \begin{equation}
        \label{eq:17}
        \int_{\Omega} \left[ \frac{1}{2} \nabla \psi^{(n+1)} \nabla v_h + \left( V(\mathbf{x}) + \beta |\psi^{(n)}|^2 + \frac{1}{\Delta t} \right) \psi^{(n+1)} v_h \right] d\mathbf{x} = \int_{\Omega} \frac{1}{\Delta t} \psi^{(n)} v_h, d\mathbf{x}, \quad \forall v_h \in V_h,
    \end{equation}
    where $V_h  = \{{v_h\in{H_0^1}}:v_h|_{\mathcal{K}}\in{P^2},\forall \mathcal{K}\in\mathcal{T}_h\}$.

    The wave function $\psi$ can be approximated as 
    \begin{equation}
        \label{eq:18}
    \psi=\sum\limits_{j=1}^{N_{basis}}{v_j}\Psi_j,
    \end{equation}
    where $N_{\text{basis}}$ denotes the dimension of the finite element space $V_h$ spanned by the basis functions, $\Psi_j$ is the coefficient associated with the $j$-th basis function, and $v_j$ represents the corresponding basis function. The basis functions are typically chosen such that $v_j$ equals one at the $j$-th finite element interpolation node and vanishes at all other nodes. Consequently, the coefficient $\Psi_j$ coincides with the value of the wave function at the corresponding node.
    Substituting Eq. \eqref{eq:18} into Eq. \eqref{eq:17} on space $V_h$, we have 
    \begin{equation}
        \label{eq:19}
        \sum_{j=1}^{n}\int_{\Omega}[\frac{1}{2}\nabla v_{j}\nabla v_{i}+(V(x)+\beta|\psi^{(n)}|^{2}+\frac{1}{\Delta t})v_{j}v_{i}]d\mathbf{x}\Psi_j^{(n+1)} = \int_{\Omega}\frac{1}{\Delta t}\psi ^{(n)} v_i d\mathbf{x}.
    \end{equation}
    As a result, the unknowns can be solved via the discretized linear equation system 
    \begin{equation}
        \label{eq:20}
        A\Psi^{(n+1)}=b,
    \end{equation}
    where
    \begin{equation}\begin{aligned}
        \label{eq:21}
        A_{ij}&=\int_{\Omega}[\frac{1}{2}\nabla v_{j}\nabla v_{i}+(V(x)+\beta|\psi^{(n)}|^{2}+\frac{1}{\Delta t})v_{j}v_{i}]d\mathbf{x}, \quad i=1,...,N,j=1,...,M,
        \\
        b_{i}&=\int_{\Omega}\frac{1}{\Delta t}\psi^{(n)} v_i d\mathbf{x}, \quad i=1,...,N.
    \end{aligned}\end{equation}

    Since the wave functions of the GPE often varies differently on the whole domain, it is necessary to use adaptive finite element methods to improve the mesh quality. In the following subsection, we will introduce the adaptive finite element method based on jump-type \textit{a posterior} error estimation \cite{Verfürth2013}.

    \subsection{Adaptive finite element method}

    Compared with the uniform mesh method, the adaptive mesh method provides higher numerical accuracy and requires fewer mesh grids. In this study, we use the $h$-adaptive method, which includes local refinement or coarsening of the mesh grids. The error indicator plays an important role in guiding the mesh refinement or coarsening. Specifically, we employ a jump-based \textit{a posterior} error estimator \cite{Verfürth2013} in which the local indicator is constructed purely from the jump of the discrete flux across the interior faces of the mesh. Let $u \in V$ and $u_{\mathcal{T}} \in V_h$ be the exact and finite element solutions of Eq.~\eqref{eq:16} and Eq.~\eqref{eq:17}, respectively. The corresponding jump-type \textit{a posteriori} error indicator for the element $\mathcal{T}_K$ is given by
    \begin{equation}
        \label{eq:22}
        \eta_{R,K}= \left(\sum_{E \in \partial \mathcal{T}_K} \frac{1}{2} h_E\sum_{l=1}^{p} \left\| \mathbb{J}_E(\psi_l) \right\|_E^2 \right)^{\frac{1}{2}},
    \end{equation}
    where $h_E$ stands for the largest length of the common face $E$ of $\mathcal{T}_K$, $\mathbb{J}_E$ is the jump term on the element $\mathcal{T}_K$ and $\mathcal{T}_J$, whose formulation is written as 
    $$
    \mathbb{J}_E(u_\mathcal{T}) = (\nabla u_l \Big|_{\mathcal{T}_K} - \nabla u_l \Big|_{\mathcal{T}_J}) \cdot \mathbf{n}_E,
    $$
    where $\mathbf{n}_E$ denotes the unit exterior normal to the element $E$.

    With the error indicator Eq. \eqref{eq:21}, the mesh adaptive algorithm is demonstrated in Algorithm \ref{alg:adaptive_mesh}.

    \begin{algorithm}[h]
            \caption{Adaptive method for solving the GPE}
            \label{alg:adaptive_mesh}
            \begin{algorithmic}[1]
                \State \textbf{Input:} Initial mesh $\mathcal{T}^{(0)}$, the tolerance between two adjacent solutions $tol_1$, adaptive tolerance $tol_2$, initial wave function $\psi^{(0)}$ satisfying the Thomas-Fermi approximation, start time $t^{(0)}=0$, initial global error estimator $\eta$ and time step $\Delta t$, $n=0$, $E_{old}= 0, E =1$. 
                \While{$|E-E_{old}|> tol_2$}
                \State $n=0, ~E_{old} = E$.
                    \While{$\|\psi^{(n+1)} - \psi^{(n)}\| > tol_1$}
                        \State \( t^{(n+1)} = t^{(n)} + \Delta t \);
                        \State Assemble the stiffness matrix $A$ and right-hand side vector $b$ on the current mesh \( \mathcal{T}^{(n)} \);
                        \State Solve the linear system $A\Psi^{(n+1)} = b$ to obtain the numerical solution $\psi^{(n+1)}$;    
                        \State Evaluate the energy $E$ using Eq. \eqref{eq:5};
                        \State $n = n+1$;      
                    \EndWhile
                     \State Compute the \textit{a posteriori} error indicator $\eta_{R,K}$ for each element $K \in \mathcal{T}^{(n)}$;                            
                    \State Mesh adaption based on the error indicator $\eta_{R,K}$ \cite{Li2005};
                    \State Update $\phi^{(n+1)}$ from the old mesh to the refined mesh;
                \EndWhile
                \State Compute the final ground state energy $E$ and chemical potential $\mu$ on the adapted mesh;
                \State \textbf{Output:} Adapted mesh \( \mathcal{T}^{*} \), ground state energy $E$, and chemical potential $\mu$.
            \end{algorithmic}
    \end{algorithm}

    \section{The multi-mesh adaptive finite element method}

    In many applications of the GPE, such as multi-component GPEs or the computation of both ground and excited states, the resulting wave functions may exhibit significantly different spatial regularities. Certain components may vary sharply in localized regions, while others remain smooth over most of the domain. In such situations, employing a single mesh to simultaneously resolve all components typically forces global refinement dictated by the most oscillatory solution, leading to excessive degrees of freedom and increased computational expense.

    To overcome this difficulty, a multi-mesh adaptive finite element method is employed \cite{Voigt2012,Kuang2024}. In this framework, each component of the system is solved on its own adaptively refined mesh, allowing the computational grid to concentrate resolution only where the corresponding wave function of the component requires it. This enables accurate representation of the distinct solution features while reducing unnecessary refinement in smooth regions, thereby significantly improving overall efficiency.

    \subsection{The model equation}
    To demonstrate this method, we consider the coupled Gross-Pitaevskii system as an example. This system consists of two GPEs:
    \begin{equation}\left\{\begin{aligned}
        \label{eq:24}
        & \mathrm{i}\partial_{t}\psi_{i}=\left[-\frac{1}{2}\nabla^2+U_i+\sum_{k=1}^2V_{ik}|\psi_{k}|^2\right]\psi_{i}, \\
        & \|\psi_{i}\|_{L^2}^2=N_i,\quad i=1,2,
    \end{aligned}\right.\end{equation}
    with scaled and off-centred harmonic potentials
    \begin{equation}
        \label{eq:25}
        U_i=\frac{1}{2}[(\omega_{1i}(x-x_i))^2+(\omega_{2i}(y-y_i))^2].
    \end{equation}

    \subsection{Discretization}

    We apply the imaginary time propagation method on Eq. \eqref{eq:24} to obtain 
    \begin{equation}\left.\left\{\begin{aligned}
        \label{eq:26}
        \partial_{t}\psi_{1}=\left[\frac{1}{2}\nabla^2-U_1 - V_{11}|\psi_{1}|^2 - V_{12}|\psi_{2}|^2\right]\psi_{1},
        \\
        \partial_{t}\psi_{2}=\left[\frac{1}{2}\nabla^2-U_2 - V_{21}|\psi_{1}|^2 - V_{22}|\psi_{2}|^2\right]\psi_{2}.
    \end{aligned}\right.\right.\end{equation}
    Next, we adopt the backward Euler scheme for temporal discretization and we have 
    \begin{equation}\left.\left\{\begin{aligned}
        \label{eq:27}
        \frac{\psi_{2}^{(n+1)} - \psi_{1}^{(n)}}{\Delta t}=\left[\frac{1}{2}\nabla^2-U_1 - V_{11}|\psi^{(n)}_{1}|^2 - V_{12}|\psi^{(n)}_{2}|^2\right]\psi_{1}^{(n+1)},
        \\
        \frac{\psi_{2}^{(n+1)} - \psi_{1}^{(n)}}{\Delta t}=\left[\frac{1}{2}\nabla^2-U_2 - V_{21}|\psi^{(n)}_{1}|^2 - V_{22}|\psi^{(n)}_{2}|^2\right]\psi_{2}^{(n+1)}.
    \end{aligned}\right.\right.\end{equation}
    Rearranging the terms we obtain 
    \begin{equation}\left.\left\{\begin{aligned}
        \label{eq:28}
        \left[-\frac{1}{2}\nabla^2 + U_1 + V_{11}|\psi^{(n)}_{1}|^2 + V_{12}|\psi^{(n)}_{2}|^2+\frac{1}{\Delta t}\right]\psi_{1}^{(n+1)}=\frac{1}{\Delta t}\psi_{1}^{(n)},
        \\
        \left[-\frac{1}{2}\nabla^2 + U_2 + V_{21}|\psi^{(n)}_{1}|^2 + V_{22}|\psi^{(n)}_{2}|^2+\frac{1}{\Delta t}\right]\psi_{2}^{(n+1)}=\frac{1}{\Delta t}\psi_{2}^{(n)}.
    \end{aligned}\right.\right.\end{equation}
    The variational formulation of this system now becomes: find $(\psi^{(1)},\psi^{(2)}) \in H_{0}^{1}(\Omega) \times H^{1}_0(\Omega)$ such that
    \begin{equation}\left.\left\{\begin{aligned}
        \label{eq:29}
        \int_{\Omega}\left[-\frac{1}{2}\nabla^2 + U_1 + V_{11}|\psi^{(n)}_{1}|^2 + V_{12}|\psi^{(n)}_{2}|^2+\frac{1}{\Delta t}\right]\psi_{1}^{(n+1)}v_{1}d\mathbf{x}=\int_{\Omega}\frac{1}{\Delta t}\psi_{1}^{(n)}v_{1}d\mathbf{x},
        \\
        \int_{\Omega}\left[-\frac{1}{2}\nabla^2 + U_2 + V_{21}|\psi^{(n)}_{1}|^2 + V_{22}|\psi^{(n)}_{2}|^2+\frac{1}{\Delta t}\right]\psi_{2}^{(n+1)}v_{2}d\mathbf{x}=\int_{\Omega}\frac{1}{\Delta t}\psi_{2}^{(n)}v_{2}d\mathbf{x}.
    \end{aligned}\right.\right.\end{equation}

    In order to discretize the equations, $\mathcal{T}_{1}$ and $\mathcal{T}_{2}$ are taken as distinct simplex partitions of the domain $\Omega$. Define
\[
V_{1,h} = \{v_{1,h} \in H_0^1 : v_{1,h}|_\mathcal{K} \in P^2, \ \forall \mathcal{K} \in \mathcal{T}_1\}, \quad
V_{2,h} = \{v_{2,h} \in H_0^1 : v_{2,h}|_\mathcal{K} \in P^2, \ \forall \mathcal{K} \in \mathcal{T}_2\}.
\]
These finite element spaces are formed by globally continuous, piecewise polynomial functions of an arbitrary fixed degree. The discrete problem is then formulated as: find $(\psi_1, \psi_2) \in V_{1,h} \times V_{2,h}$ such that

    \begin{footnotesize}\begin{equation}\left.\left\{\begin{aligned}
        \label{eq:30}
        \int_{\Omega}\left[\frac{1}{2}\nabla\psi_{1}^{(n+1)} \nabla{v_{1,h}} + (U_1 + V_{11}|\psi^{(n)}_{1}|^2 + V_{12}|\psi^{(n)}_{2}|^2+\frac{1}{\Delta t})\psi_{1}^{(n+1)}v_{1,h}\right]d\mathbf{x}=\int_{\Omega}\frac{1}{\Delta t}\psi_{1}^{(n)}{v_{1,h}}d\mathbf{x},\quad \forall v_{h}^{1} \in V_{h}^{1}(\Omega),
        \\
        \int_{\Omega}\left[\frac{1}{2}\nabla\psi_{2}^{(n+1)} \nabla{v_{2,h}} + (U_2 + V_{21}|\psi^{(n)}_{1}|^2 + V_{22}|\psi^{(n)}_{2}|^2+\frac{1}{\Delta t})\psi_{2}^{(n+1)}v_{2,h}\right]d\mathbf{x}=\int_{\Omega}\frac{1}{\Delta t}\psi_{2}^{(n)}{v_{2,h}}d\mathbf{x},\quad \forall v_{h}^{2} \in V_{h}^{2}(\Omega).
    \end{aligned}\right.\right.\end{equation}\end{footnotesize}

    In terms of spatial discretization, let us define \{$v_{1,i}|1 \leq i \leq N$\} and \{$v_{2,j}|1 \leq j \leq M$\} to be the nodal basis of $V_{1,h}$ and $V_{2,h}$, respectively. Thus the wave function $\psi_j$ can be approximated as 
    \begin{equation}
        \label{eq:31}
        \psi^{(n+1)}_{1,j}=\sum_{j=1}^{N}v_{1,j}\Psi^{(n+1)}_{1,j},
    \end{equation}
    \begin{equation}
        \label{eq:32}
        \psi^{(n+1)}_{2,j}=\sum_{j=1}^{M}v_{2,j}\Psi^{(n+1)}_{2,j}.
    \end{equation}

    To find an approximate solution, using the above representations, Eq. \eqref{eq:30} can be rewritten as  
    \begin{footnotesize}\begin{equation}\left.\left\{\begin{aligned}
        \label{eq:33}
        \sum_{j=1}^{N} \{\int_{\Omega}\left[\frac{1}{2}\nabla{v_{1,j}} \nabla{v_{1,i}} + (U_1 + V_{11}|\psi^{(n)}_{1}|^2 + V_{12}|\psi^{(n)}_{2}|^2 + \frac{1}{\Delta t}){v_{1,j}}{v_{1,i}} \right]d\mathbf{x} \Psi_{1,j}^{(n+1)} 
        = \int_{\Omega}\frac{1}{\Delta t}\psi_{1}^{(n)}{v_{1,i}}d\mathbf{x},
        \\
        \sum_{j=1}^{M} \{\int_{\Omega}\left[\frac{1}{2}\nabla{v_{2,j}} \nabla{v_{2,i}} + (U_2 + V_{21}|\psi^{(n)}_{1}|^2 + V_{22}|\psi^{(n)}_{2}|^2 + \frac{1}{\Delta t}){v_{2,j}}{v_{2,i}} \right]d\mathbf{x} \Psi_{2,j}^{(n+1)} 
        = \int_{\Omega}\frac{1}{\Delta t}\psi_{2}^{(n)}{v_{2,i}}d\mathbf{x}.
    \end{aligned}\right.\right.\end{equation}\end{footnotesize}
    The above equations lead to the following linear systems: 
    \begin{equation}\left.\left\{\begin{aligned}
        \label{eq:34}
        A_1\Psi_{1}^{(n+1)}={b}^1,
        \\
        A_2\Psi_{2}^{(n+1)}={b}^2,
    \end{aligned}\right.\right.\end{equation}
    where
    \begin{equation}\left.\left.\left.\left\{\begin{aligned}
        \label{eq:35}
        A_{1,ij}&=\int_{\Omega}\left[\frac{1}{2}\nabla{v_{1,j}} \nabla{v_{1,i}} + (U_1 + V_{11}|\psi^{(n)}_{1}|^2 + V_{12}|\psi^{(n)}_{2}|^2 + \frac{1}{\Delta t}){v_{1,j}}{v_{1,i}}\right]d\mathbf{x}, 
        \\
        A_{2,ij}&=\int_{\Omega}\left[\frac{1}{2}\nabla{v_{2,j}} \nabla{v_{2,i}} + (U_2 + V_{21}|\psi^{(n)}_{1}|^2 + V_{22}|\psi^{(n)}_{2}|^2+\frac{1}{\Delta t}){v_{2,j}}{v_{2,i}}\right]d\mathbf{x}, 
        \\
        b_{1,i}&=\int_{\Omega}\frac{1}{\Delta t}\psi_{1}^{(n)}{v_{1,i}}d\mathbf{x},
        \\
        b_{2,i}&=\int_{\Omega}\frac{1}{\Delta t}\psi_{2}^{(n)}{v_{2,i}}d\mathbf{x}.
    \end{aligned}\right.\right.\right.\right.\end{equation}
    The energy functional of the coupled system is expressed as follows \cite{Caliari2008}
    $$E(\psi_1, \psi_2) = E_{\infty}(\psi_1, \psi_2) + V_{12} \int_{\Omega} |\psi_1|^2 |\psi_2|^2d\mathbf{x},$$
    where
    $$E_{\infty}(\psi_1, \psi_2) = \sum_{i=1}^{2} E_i(\psi_i(\mathbf{x},t)),$$
    and 
    $$E_i(\psi_i(\mathbf{x},t)) = \frac{1}{2} \int_{\Omega} |\nabla \psi_i(\mathbf{x}, t)|^2d\mathbf{x} + \int_{\Omega} U_i(\mathbf{x}) |\psi_i(\mathbf{x}, t)|^2d\mathbf{x} + \frac{V_{ii}}{2} \int_{\Omega} |\psi_i(\mathbf{x}, t)|^4d\mathbf{x}, \quad i=1,2. $$

    \subsection{The hierarchical geometry tree}
    An efficient management mechanism requires a well-designed data structure for mesh grids. Our algorithm employs the hierarchy geometry tree (HGT) \cite{Li2005,Bao2012,Cai2024,Kuang2024}. In this hierarchical representation, all geometries exhibit belonging-to relationships. For instance, if a triangular face forms part of a tetrahedron, it belongs to that tetrahedron. Likewise, that triangle possesses its own edges and vertices. This layered organization allows flexible referencing of geometric information for triangles and supports efficient mesh refinement and coarsening operations.

    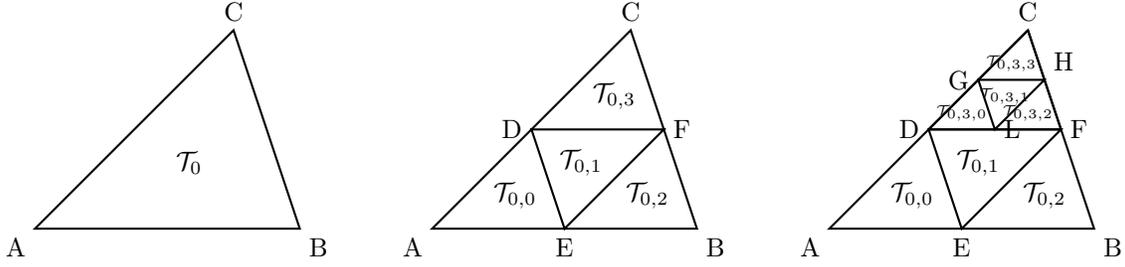
\begin{figure}[htbp]
    \centering
        \begin{tikzpicture}[scale=0.88, baseline]

            \begin{scope}[shift={(0,0)}]
                \coordinate (A) at (-3,0); \coordinate (B) at (1,0); \coordinate (C) at (0,3);
                \draw[thick] (A)--(B)--(C)--cycle;
                \node[below left]  at (A) {A}; \node[below right] at (B) {B}; \node[above] at (C) {C};
                \node at (-0.67,1) {$\mathcal{T}_0$};
            \end{scope}

            \begin{scope}[shift={(6,0)}]
                \coordinate (A) at (-3,0); \coordinate (B) at (1,0); \coordinate (C) at (0,3);
                \coordinate (D) at ($ (A)!0.5!(C) $); 
                \coordinate (E) at ($ (A)!0.5!(B) $); 
                \coordinate (F) at ($ (B)!0.5!(C) $); 
                \draw[thick] (A)--(B)--(C)--cycle (D)--(E)--(F)--cycle;
                \node[below left]  at (A) {A}; \node[below right] at (B) {B}; \node[above] at (C) {C};
                \node[left]  at (D) {D}; \node[below] at (E) {E}; \node[right] at (F) {F};
                \node at (-1.75,0.5) {$\mathcal{T}_{0,0}$};
                \node at (-0.75,1.0) {$\mathcal{T}_{0,1}$};
                \node at (0.25,0.5) {$\mathcal{T}_{0,2}$};
                \node at (-0.25,2.0) {$\mathcal{T}_{0,3}$};
            \end{scope}

            \begin{scope}[shift={(12,0)}]
                \coordinate (A) at (-3,0); \coordinate (B) at (1,0); \coordinate (C) at (0,3);
                \coordinate (D) at ($ (A)!0.5!(C) $); 
                \coordinate (E) at ($ (A)!0.5!(B) $); 
                \coordinate (F) at ($ (B)!0.5!(C) $); 

                \coordinate (G) at ($ (D)!0.5!(C) $); 
                \coordinate (H) at ($ (C)!0.5!(F) $); 
                \coordinate (L) at ($ (D)!0.5!(F) $); 

                \draw[thick] (A)--(B)--(C)--cycle;
                \draw[thick] (D)--(E)--(F)--cycle;

                \draw[thick] (G)--(H)--(L)--cycle;
                \draw[thick] (D)--(G) (C)--(G) (C)--(H) (F)--(H) (D)--(L) (F)--(L);

                \node[below left]  at (A) {A}; \node[below right] at (B) {B}; \node[above] at (C) {C};
                \node[left]  at (D) {D}; \node[below] at (E) {E}; \node[right] at (F) {F};
                \node[left]  at (G) {G}; \node[above right] at (H) {H}; \node[right] at (L) {L};

                \node[font=\tiny] at (-1.0,1.75) {$\mathcal{T}_{0,3,0}$};
                \node[font=\tiny] at (-0.35,2.0) {$\mathcal{T}_{0,3,1}$};
                \node[font=\tiny] at (0,1.75) {$\mathcal{T}_{0,3,2}$};
                \node[font=\tiny] at (-0.25,2.5) {$\mathcal{T}_{0,3,3}$};

                \node at (-1.75,0.5) {$\mathcal{T}_{0,0}$};
                \node at (-0.75,1.0) {$\mathcal{T}_{0,1}$};
                \node at (0.25,0.5) {$\mathcal{T}_{0,2}$};
            \end{scope}

        \end{tikzpicture}
    \caption{The refinement of a triangle $\mathcal{T}_0$ and its sub-triangle $\mathcal{T}_{0,3}$.}
    \label{fig:hgt}
    \end{figure}

\begin{figure}[htbp]
    \centering
        \begin{tikzpicture}[scale=0.8, transform shape, 
        level distance=1.2cm,
        level 1/.style={sibling distance=2.5cm},
        level 2/.style={sibling distance=1.2cm}]
  
        \node {$\mathcal{T}_{0}$}
        child {node {$\mathcal{T}_{0,0}$}}
        child {node {$\mathcal{T}_{0,1}$}}
        child {node {$\mathcal{T}_{0,2}$}}
        child {node {$\mathcal{T}_{0,3}$}
        child {node {$\mathcal{T}_{0,3,0}$}}
        child {node {$\mathcal{T}_{0,3,1}$}}
        child {node {$\mathcal{T}_{0,3,2}$}}
        child {node {$\mathcal{T}_{0,3,3}$}}
        };

        \end{tikzpicture}
    \caption{The quadtree data structure for the mesh in the right of Figure \ref{fig:hgt}.}  
    \label{Hierarchical geometry tree}
    \end{figure}
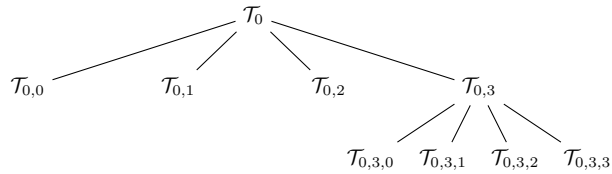

    Let us examine a regular triangulation $\mathcal{T}_0$ over a two-dimensional domain $\Omega$. A triangle within this triangulation can be subdivided into four smaller triangles as depicted in Figure \ref{fig:hgt}. Subsequently, by further refining $\mathcal{T}_{0,3}$, an additional set of sub-triangles $\{\mathcal{T}_{0,3,0},\mathcal{T}_{0,3,1},\mathcal{T}_{0,3,2},\mathcal{T}_{0,3,3}\}$ is created, as illustrated in Figure \ref{fig:hgt}. Repeating this procedure iteratively results in a sequence of triangulations, forming a four-branch tree structure rooted at each initial triangle, as shown in Figure \ref{Hierarchical geometry tree}.

\begin{figure}[htbp]
    \centering
        \begin{tikzpicture}[scale=0.88, baseline]
            \begin{scope}[shift={(0,0)}]
                \coordinate (D) at (-3,0);
                \coordinate (F) at (1,0);
                \coordinate (E) at (0,3);
                \coordinate (L) at ($ (D)!0.5!(F) $);  

                \draw[thick] (D)--(E)--(F)--cycle;
                \draw[dashed] (E)--(L);                

                \node[below left]  at (D) {D};
                \node[below right] at (F) {F};
                \node[above]       at (E) {E};
                \node[below]       at (L) {L};
            \end{scope}
        \end{tikzpicture}
        \caption{Twin-triangle.}
        \label{fig:twin-triangle}
    \end{figure}
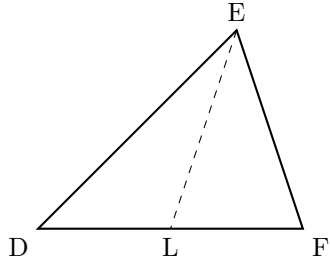

    Efficient mesh management is achieved in the HGT via a hierarchical geometric representation and an associated tree data structure. However, directly constructing finite element space on the mesh will generate hanging points on the adjacent edges of the refined triangles. If a refined element has multiple hanging points, it will be global refined again. If there is only one hanging point, for example, in the refinement process of $\triangle DFC$ in Figure \ref{fig:hgt}, the hanging point $L$ appears. For dealing with these kind of hanging points, the twin-triangle geometry is introduced, which is demonstrated in Figure \ref{fig:twin-triangle}.

    There are precisely three possible relationships between any two elements: equality, belonging-to, or no overlap. Interactions among elements that are equal or non-overlap is straightforward. Therefore, our focus must be directed toward the second type of relationship. In the present method, the information of the two wave functions in the coupled Gross-Pitaevskii system come from two different meshes $\mathcal{T}_1$ and $\mathcal{T}_2$, and their treatment primarily revolves around numerical integration. To illustrate, we consider the computation of the coupling term presented in Eq. (\ref{eq:35}). As example, it can be written as
    \begin{equation}
        \label{eq:36}
        \int_{\Omega} V_{12}|\psi^{(n)}_{2}|^2 v_{1,j} v_{1,i}d\mathbf{x} .
    \end{equation}

    \begin{figure}[htbp]
    \centering
        \begin{tikzpicture}[scale=0.88, baseline]
            \begin{scope}[shift={(0,0)}]
                \coordinate (A) at (-3,0); \coordinate (B) at (1,0); \coordinate (C) at (0,3);
                \coordinate (D) at (-1.5,1.5); \coordinate (E) at (-1,0);   \coordinate (F) at (0.5,1.5);
                \coordinate (G) at (-2.25,0.75); \coordinate (H) at (-2,0); \coordinate (I) at (-1.25,0.75);
                \draw[dashed] (I) -- (F);
                \draw[thick] (A)--(B)--(C)--cycle
                             (D)--(E)--(F)--cycle
                             (G)--(H)--(I)--cycle;
                \node[below left]  at (A) {A}; \node[below right] at (B) {B}; \node[above] at (C) {C};
                \node[left]  at (D) {D}; \node[below] at (E) {E}; \node[right] at (F) {F};
                \node[left]  at (G) {G}; \node[below] at (H) {H}; \node[right] at (I) {I};
            \end{scope}

            \begin{scope}[shift={(8,0)}]
                \coordinate (A) at (-3,0); \coordinate (B) at (1,0); \coordinate (C) at (0,3);
                \coordinate (D) at (-1.5,1.5); \coordinate (E) at (-1,0); \coordinate (F) at (0.5,1.5);

                \coordinate (J) at (-0.75,2.25); 
                \coordinate (K) at (0.25,2.25);  
                \coordinate (L) at (-0.5,1.5);  

                \draw[dashed] (E) -- (L);
                \draw[thick] (A)--(B)--(C)--cycle
                             (D)--(E)--(F)--cycle
                             (J)--(K)--(L)--cycle;
                \node[below left]  at (A) {A}; \node[below right] at (B) {B}; \node[above] at (C) {C};
                \node[left]  at (D) {D}; \node[below] at (E) {E}; \node[right] at (F) {F};
                \node[left]  at (J) {J}; \node[right] at (K) {K}; \node[below right] at (L) {L};
            \end{scope}

        \end{tikzpicture}
    \caption{Hierarchical geometry tree: $\mathcal{T}_1$(left),$\mathcal{T}_2$(right).}
    \label{fig:h}
    \end{figure}
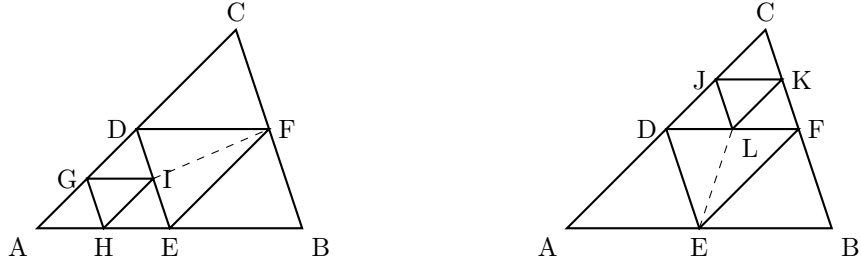

    The integration must be handled with caution, because the wave function $\psi_2^{(n)}$ and basis function $v_{1,j}, v_{1,i}$ are established on different finite element spaces. Here is an intuitive example to illustrate the proposed method for calculating Eq. \eqref{eq:36}. Assuming that the finite element space of wave function $\psi_1$ is established on mesh $\mathcal{T}_1$, and the finite element space of wave function $\psi_2$ is established on mesh $\mathcal{T}_2$, since $v_{1,j}$, $v_{1,i}$ is the basis function of wave function $\psi_1$, see Figure \ref{fig:h}. For the element $\triangle BEF$ that belongs to both $\mathcal{T}_1$ and $\mathcal{T}_2$, the integral can be calculated directly. For the remaining elements, for example, triangle $\triangle AED$ is refined in $\mathcal{T}_1$, while in $\mathcal{T}_2$, triangle $\triangle AED$ is preserved, and the same applies to triangle $\triangle CDF$.

    To ensure accurate numerical results, we employ a strategy that fully exploits the available integration points in the numerical quadrature \cite{Kuang2024}. In practice, the integral over $\triangle AED$ is evaluated through a decomposition into four refined subtriangles,
    \begin{align*}
        & V_{12}\int_{\triangle AED}  |\psi^{(n)}_{2}|^2 v_{1,j} v_{1,i} \, d\mathbf{x} \\
        &= V_{12} \sum_{k=1}^4 \int_{\triangle k} |\psi^{(n)}_{2}|^2 v_{1,j} v_{1,i} \, d\mathbf{x} \\
        &\approx V_{12} \sum_{k=1}^4 \sum_{l=1}^q \text{area}({\triangle_{k}}) J_l^k w_l^k \cdot \underbrace{|\psi^{(n)}_{2}(x_l^k)|^2}_{element {\triangle k}} \cdot \underbrace{v_{1,j}(x_l^k) v_{1,i}(x_l^k)}_{element \triangle AED}
    \end{align*}
    where $k$ denotes the four subtriangles of $\triangle AED$, $q$ is the number of integration points, $\mathbf{x}_l^k$ represents the $l$-th quadrature point on subtriangle $k$, ${J}_l^k$ is the jacobian at $\mathbf{x}_l^k$, and $w_l^k$ is the  corresponding quadrature weight. With the HGT, the solution update can be implemented efficiently.

\begin{figure}[htbp]
        \centering \small
        \begin{tikzpicture}[node distance=1.2cm, scale=0.8, every node/.style={transform shape}]
            \node (start) [startstop] {Initial mesh $\mathcal{T}_1 = \mathcal{T}_2 = \mathcal{T}$};
            \node (guess) [startstop, below of=start, yshift=-0.2cm] {Initial guess $\Psi_1^{(0)},\Psi_2^{(0)}$};
            \node (calcphi) [processBoth, below of=guess ,yshift=-0.6cm] {Calculate $\psi_1^{(n)},\psi_2^{(n)}$};
            \node (constructA) [processBoth, below of=calcphi ,yshift=-0.2cm] {Construct $A_1,A_2$};
            \node (solve) [processBoth, below of=constructA ,yshift=-0.2cm] {Solve $A_{1}\Psi_1^{(n)} = b_1,A_{2}\Psi_2^{(n)} = b_2$};
            \node (update) [processBoth, below of=solve ,yshift=-0.2cm] {Update $\Psi_1^{(n+1)},\Psi_2^{(n+1)}$};
            \node (check1) [decision, below of=update, yshift=-1.5cm] { $|\psi^{(n+1)}_i - \psi^{(n)}_i|<tol_1$ i=1,2 or $t = T_{End}$ };
            \node (check2) [decision, right of=check1, xshift=5.5cm] {$|E-E_{old}|< tol_2$};
            \node (output) [startstop, below of=check2, yshift=-1.5cm] {Output};
            \node (adapt0) [process0, above of=check2, yshift=3.5cm, xshift=-1.7cm] {Adapt mesh $\mathcal{T}_1$ with $\eta_1$};
            \node (adapt1) [process1, right of=adapt0, xshift=3.5cm] {Adapt mesh $\mathcal{T}_2$ with $\eta_2$};

            \coordinate (mergeAdapt) at ($(adapt0.north)!0.5!(adapt1.north)+(0,0.8)$);
            \node (mergepoint) at ($(guess.south)!0.5!(calcphi.north)+(0,0.1)$) {};
            \draw [arrow] (start) -- (guess);
            \draw [arrow] (guess) -- (calcphi);
            \draw [arrow] (calcphi) -- (constructA);
            \draw [arrow] (constructA) -- (solve);
            \draw [arrow] (solve) -- (update);
            \draw [arrow] (update) -- (check1);
            \draw [arrow] (check1) -- node[above] {Yes} (check2);
            \draw [arrow] (check1.west) -- ++(-2,0) node[left] {No} -- ++(0,6.9) -- (calcphi.west);
            \draw [arrow] (check2.south) -- node[right] {Yes} (output);
            \draw [arrow] (check2.north) -- ++(0,0.8) node[near start,right] {No} -| (adapt0.south);
            \draw [arrow] (check2.north) -- ++(0,0.8) node[near start,right] {No} -| (adapt1.south);
            \draw [arrow] (adapt0.north) -- ++(0,0.6) -| (mergeAdapt);
            \draw [arrow] (adapt1.north) -- ++(0,0.6) -| (mergeAdapt);
            \draw [arrow] (mergeAdapt) |- ($(mergepoint)+(0,0)$);
        \end{tikzpicture}
        \caption{Flowchart of the multi-mesh adaptive algorithm for the coupled Gross-Pitaevskii system. $E_{old}$ is the previous energy, $tol_1$ is the difference between two adjacent solutions and $tol_2$ is adaptive tolerance; $\eta_1$, $\eta_2$ are the error indicators for mesh adaptation of $\mathcal{T}_1$ and $\mathcal{T}_2$, respectively.}
        \label{fig:flowchart-multimesh}
    \end{figure}
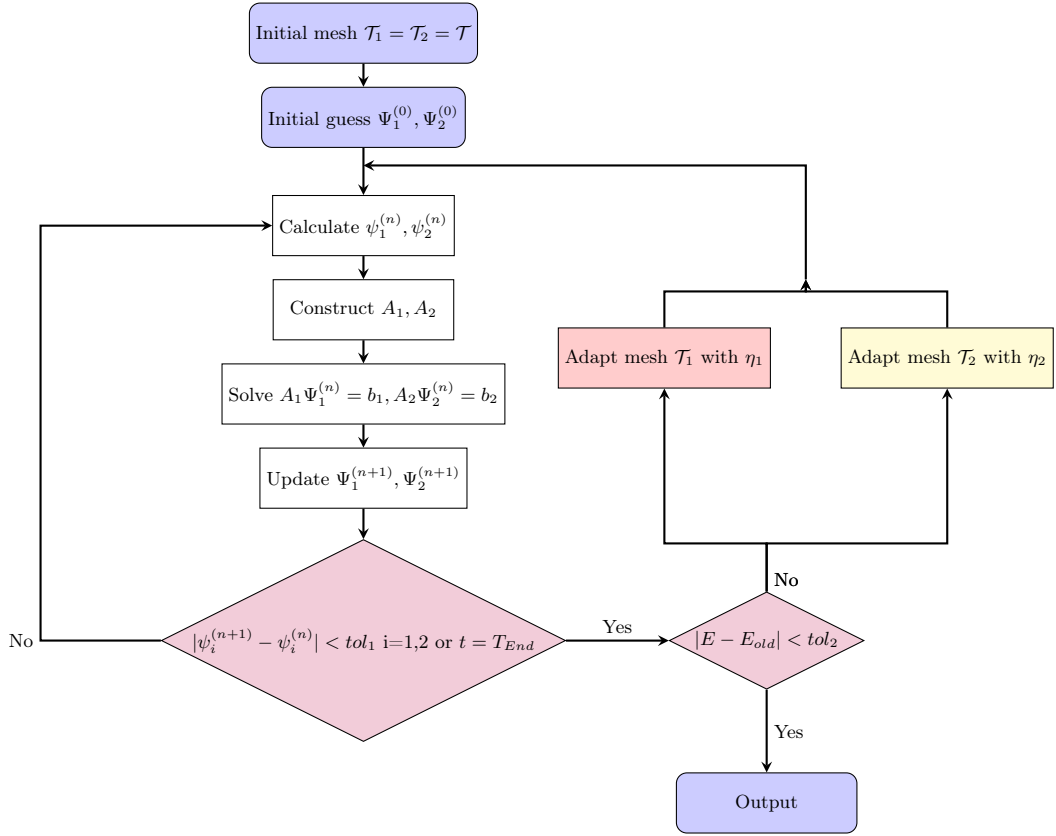

    \subsection{The multi-mesh adaptive algorithm}
    By utilizing the HGT data structure, a multi-mesh adaptive finite element method can be effectively implemented. In order to solve the coupled Gross-Pitaevskii system, we propose a multi-mesh adaptive finite element flowchart as shown in Figure \ref{fig:flowchart-multimesh}. Two meshes are adaptively adjusted using the corresponding error indicators. The mesh adaptation continues until the energy difference between two adjacent solutions is less than the prescribed tolerance.

    \section{Numerical experiments}
    In this section, the effectiveness of the proposed method is tested through three numerical experiments. In the first example, we simulated BEC and clearly observed the advantages of the adaptive finite element method. In the second example, we simultaneously calculated the ground state and first excited state of the BEC and compared it with the single-mesh method. The experimental results showed that the multi-mesh adaptive finite element method can significantly improve computational efficiency. In the final example, we used the proposed multi-mesh adaptive finite element method to solve the coupled Gross-Pitaevskii system, and the results successfully demonstrated that the proposed method can also be used to handle coupling problems. In the experiment, we all used quadratic Lagrangian element. The time iterations are terminated when the numerical solution satisfies the following criterion:
    $$\frac{\|\Psi^{(n+1)} - \Psi^{(n)}\|_{\infty}}{\Delta t} < 10^{-6}.$$

    \subsection*{Example 1}
    In the first example, we used adaptive finite element method to calculate the ground state of BEC under the combined action of harmonic and optical lattice potential
    \begin{equation}
        V(x,y) = \tfrac12(x^2 + y^2) + \kappa\!\left[\sin^2\!\left(\tfrac{\pi x}{4}\right)+\sin^2\!\left(\tfrac{\pi y}{4}\right)\right],
        \label{eq:37}
    \end{equation}
    with $\beta=1000$, $\kappa=100$ on $\Omega=[-16,16]^2$. This example is from \cite{Bao2006}. We choose $\Delta t = 0.01$. The energy and chemical potential of ground state are computed as follows: $E(\phi_g)=51.22032800,\mu(\phi_g)=66.249021.$

\begin{figure}[htbp]         
        \centering
        \setlength{\tabcolsep}{0.5pt}
        \begin{tabular}{cc}
            \includegraphics[width=0.45\linewidth]{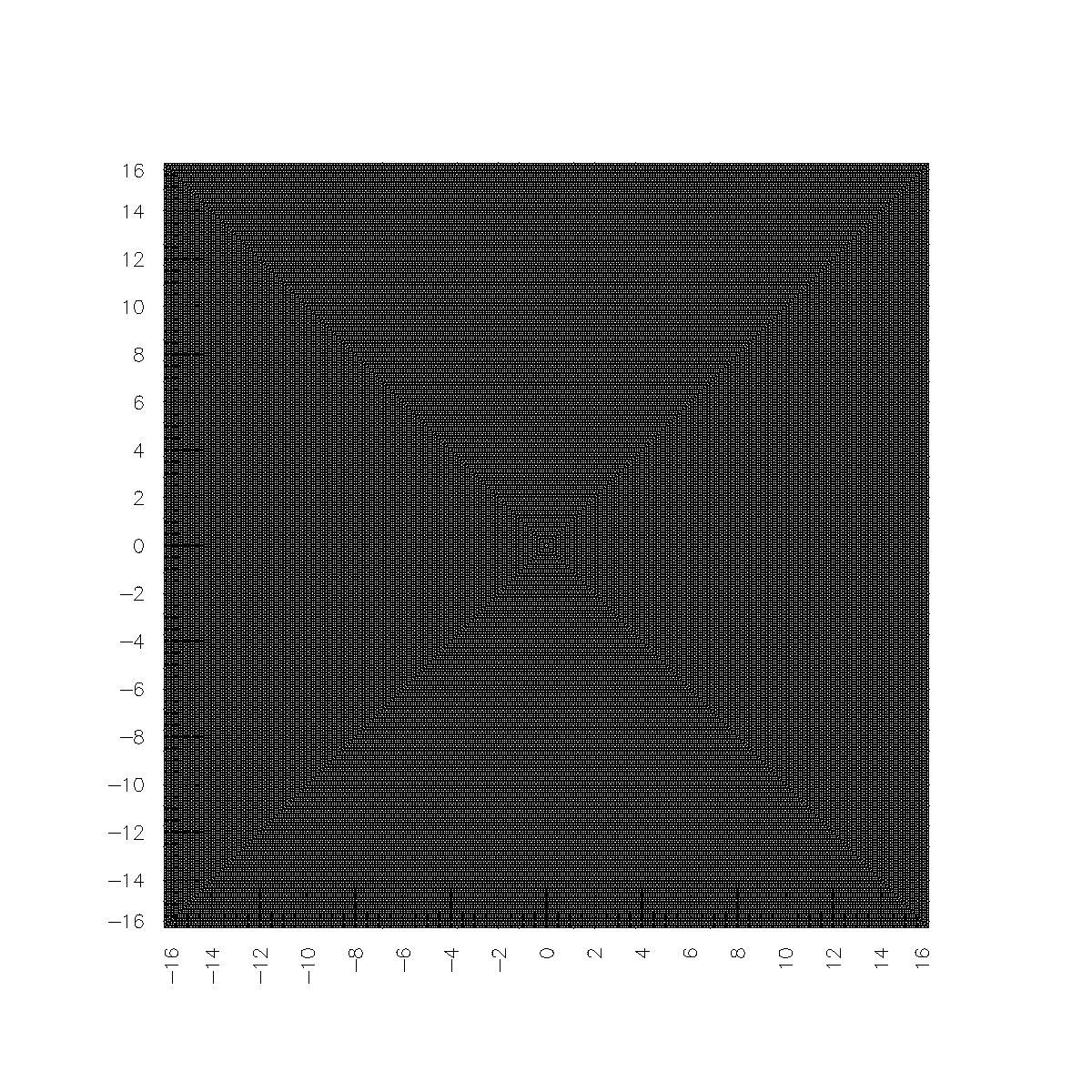}&
            \includegraphics[width=0.45\linewidth]{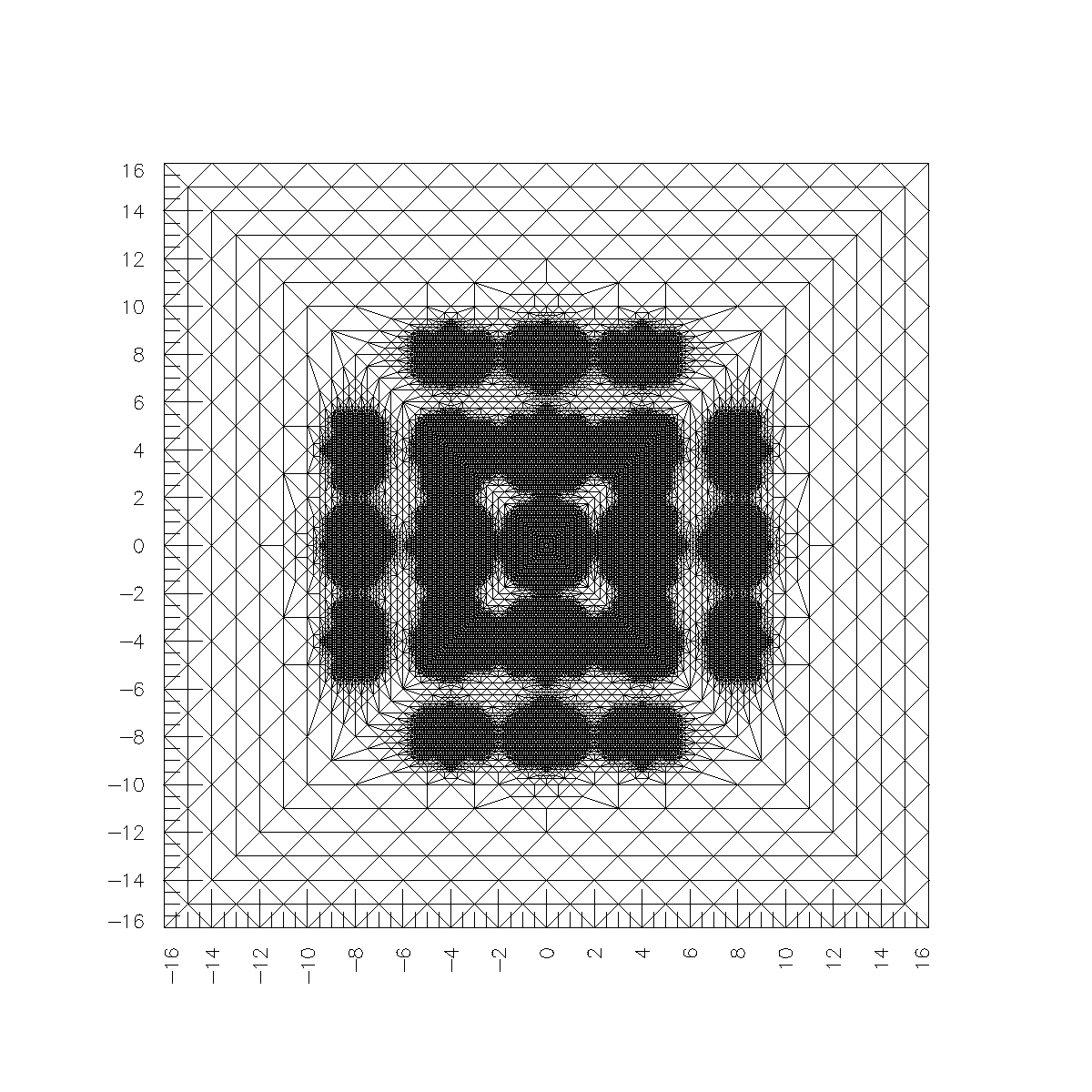}\\[-2pt]
            \small (a) uniform mesh & \small (b) adaptive mesh\\[4pt]
            \includegraphics[width=0.45\linewidth]{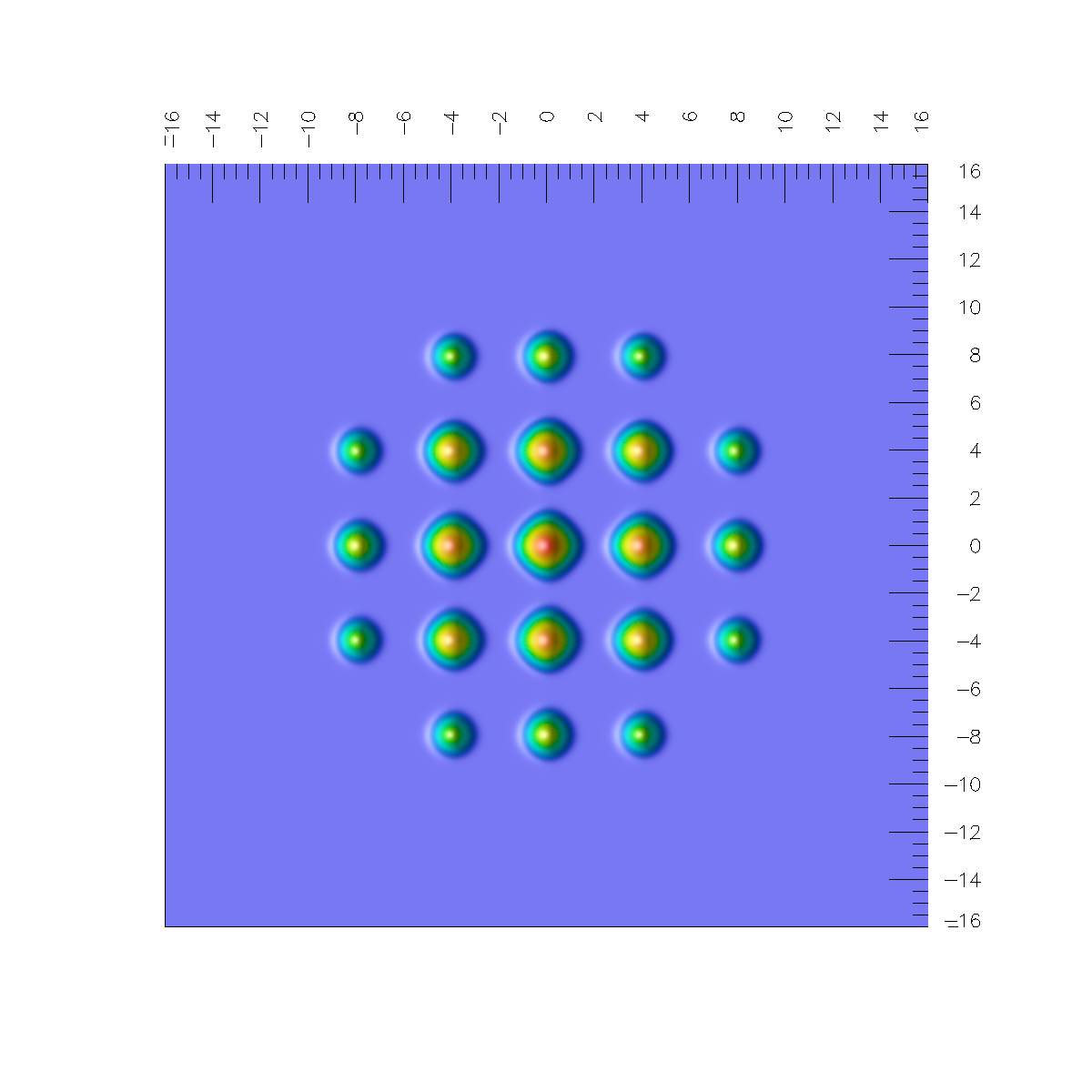}&
            \includegraphics[width=0.45\linewidth]{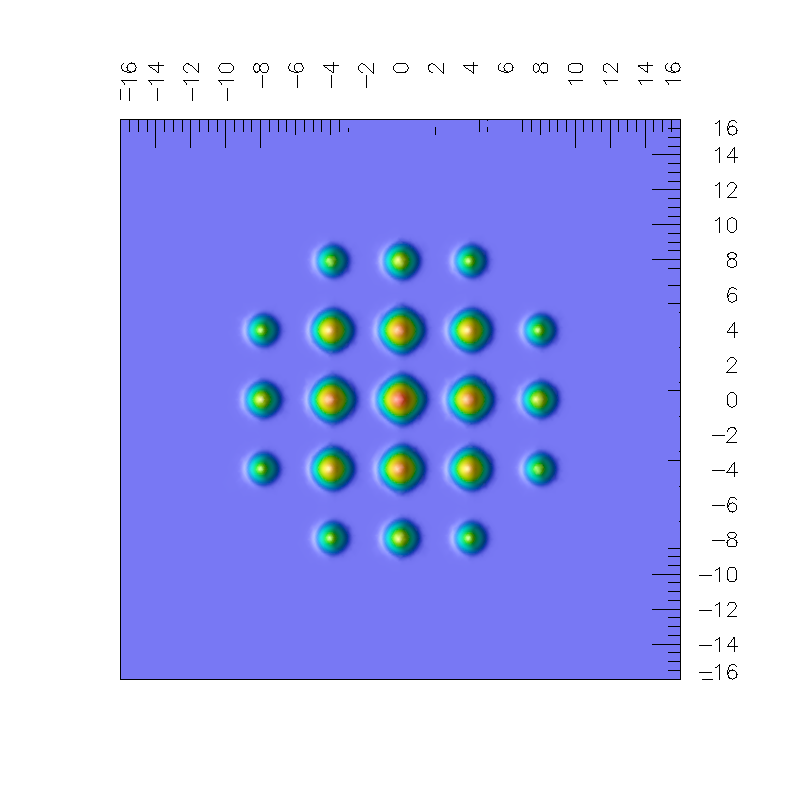}\\[-2pt]
            \small (c) uniform solution & \small (d) adaptive solution
        \end{tabular}
        \caption{Mesh and solution of uniform mesh and adaptive method in ground state for Example 1: (left: uniform mesh (525,313 dofs) and solution; right: adaptive mesh (165,889 dofs) and solution).}
        \label{adaptivegpe}
    \end{figure}

    The uniform mesh is shown in the top left of Figure \ref{adaptivegpe}, the mesh of the adaptive finite element method is shown in the top right of Figure \ref{adaptivegpe}, and the corresponding solutions are shown in the bottom of Figure \ref{adaptivegpe}.

    \begin{table}[htbp]
        \centering 
        \vspace{0.5em}
        \caption{Comparison of numerical results of uniform mesh and adaptive method in ground state, the referenced values \cite{Bao2006} are $E(\phi_g)=51.22028604,\mu(\phi_g)=66.249010$.}
        \vspace{0.5em}
        \label{tab:comparison1}
        \begin{tabular}{llllll} 
            \toprule
            {} & dofs & Energy $E$ & $\delta E$ & Chemical potential $\mu$ & $\delta \mu$ \\
            \midrule
            Uniform mesh    & 525,313 & 51.22035689 & 0.00007085 & 66.249044 & 0.000034 \\
            Adaptive method & 165,889 & 51.22032800 & 0.00004196 & 66.249021 & 0.000011 \\
            \bottomrule
        \end{tabular}
    \end{table}

    We use the numerical experiment in \cite{Bao2006} as a reference value and compare the obtained numerical results with experiments under uniform mesh. From Table \ref{tab:comparison1}, $\delta E$ and $\delta \mu$ represent the absolute errors of energy and chemical potential with respect to the reference value, respectively. Although increasing the uniform mesh to 525,313 degrees of freedom yields satisfactory accuracy, this requirement substantially exceeds the 165,889 degrees of freedom employed by the adaptive method. As shown in Table \ref{tab:comparison1}, these results underscore the efficiency of the proposed adaptive finite element method, which achieves the target numerical accuracy with markedly fewer degrees of freedom.

    \subsection*{Example 2} 
    In this example, we used the multi-mesh adaptive finite element method to simultaneously the ground state and first excited state of BEC. We use the initial values introduced in Section \ref{subsec:2.1} for solving the ground state, and also use the combined action of harmonic and optical lattice potential Eq. \eqref{eq:37}. We take $\frac{\sqrt{2}x}{\pi^{1/2}} e^{-(x^2+y^2)/2}$, $\frac{\sqrt{2}y}{\pi^{1/2}} e^{-(x^2+y^2)/2}$ for the first excited state in $x$-direction $\phi_{10}$, $y$-direction $\phi_{01}$ as the initial values \cite{Bao2006,Yang2025}, respectively. The problems are solved with $\beta = 3200, \kappa = 50$ on $\Omega=[-16,16]^2$. This example is from \cite{Zhang2016}. 
    \begin{figure}[htbp]
        \centering
        \includegraphics[width=0.5\linewidth]{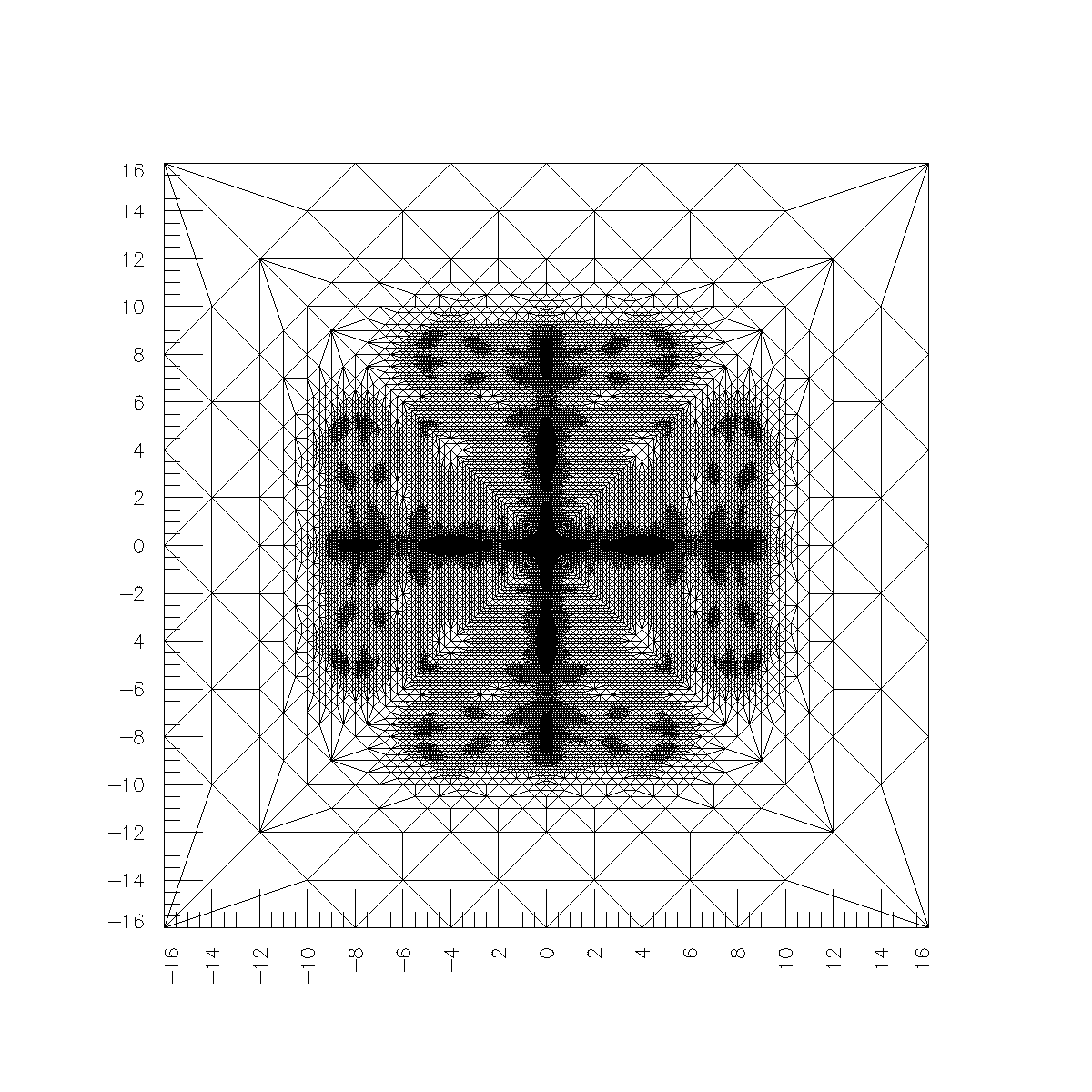}
        \caption{Simultaneously solving the desired states of a single-mesh in Example 2.}
        \label{single_excited_mesh}
    \end{figure}
    \begin{figure}[htbp]
        \centering
        \includegraphics[width=0.3\linewidth]{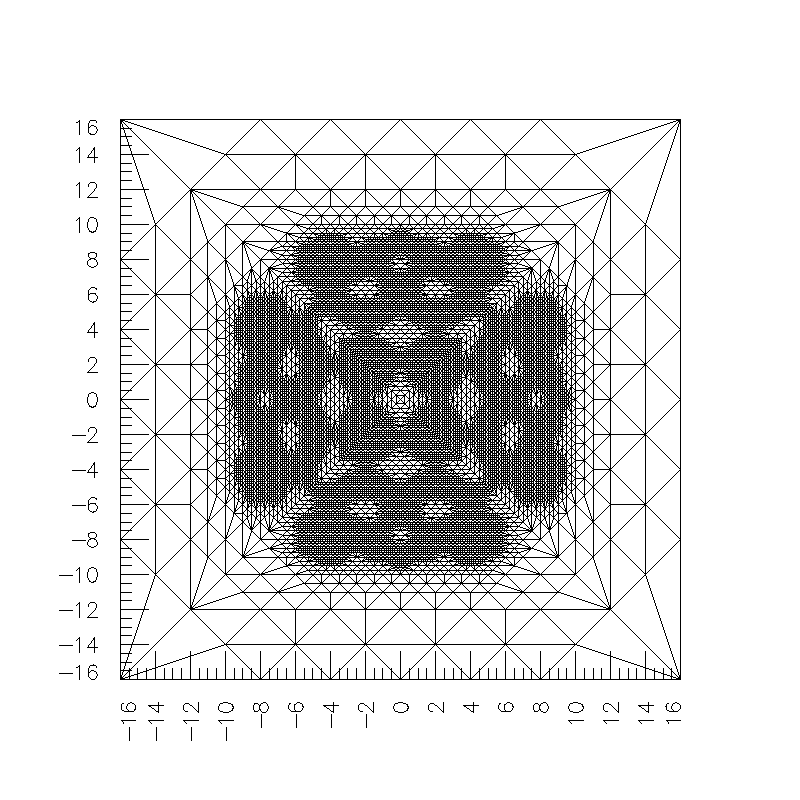}
        \includegraphics[width=0.3\linewidth]{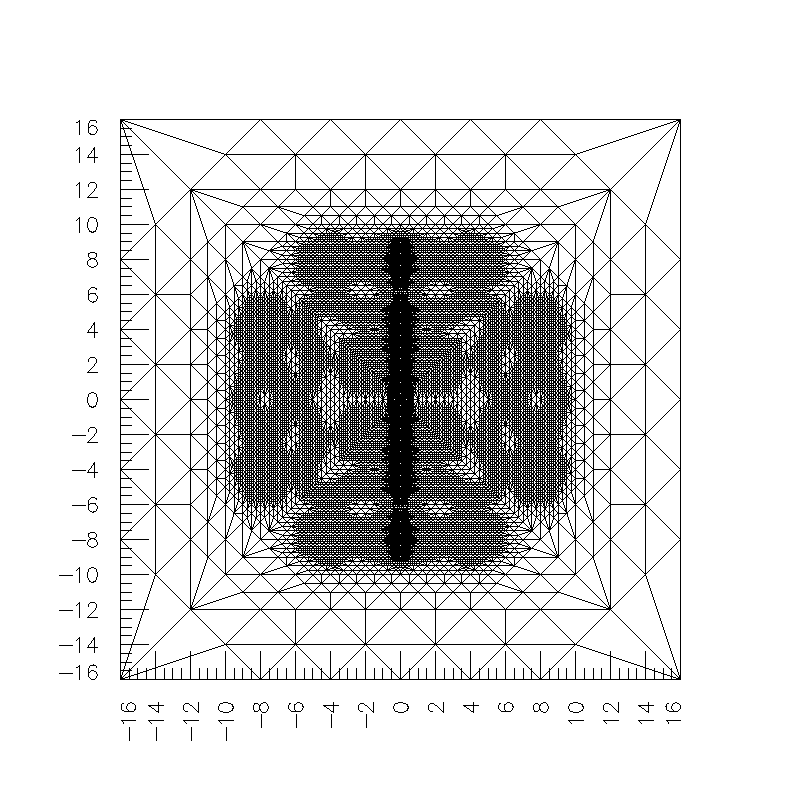}
        \includegraphics[width=0.3\linewidth]{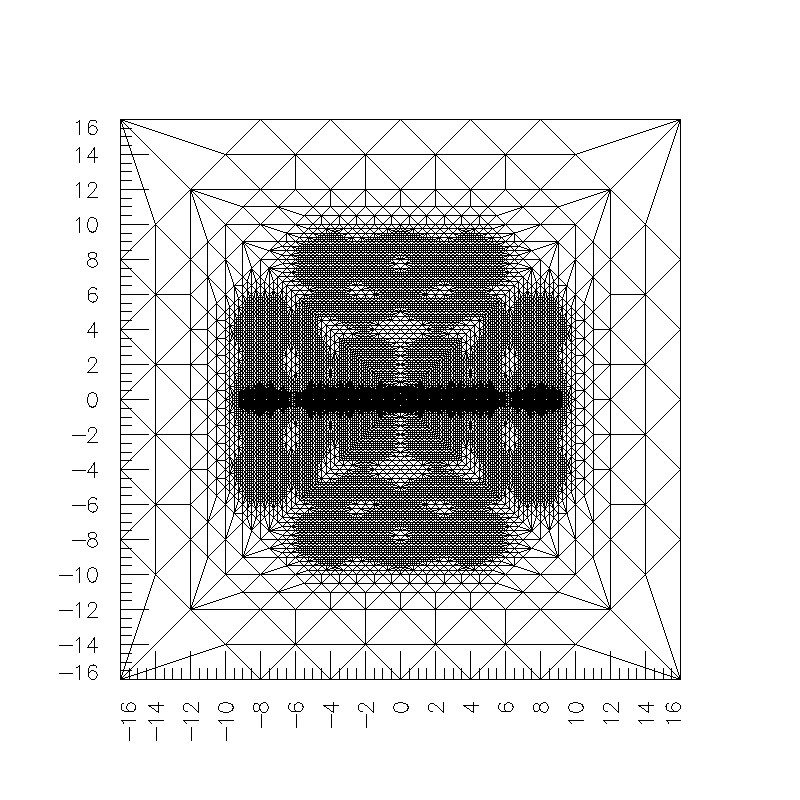}
        \includegraphics[width=0.3\linewidth]{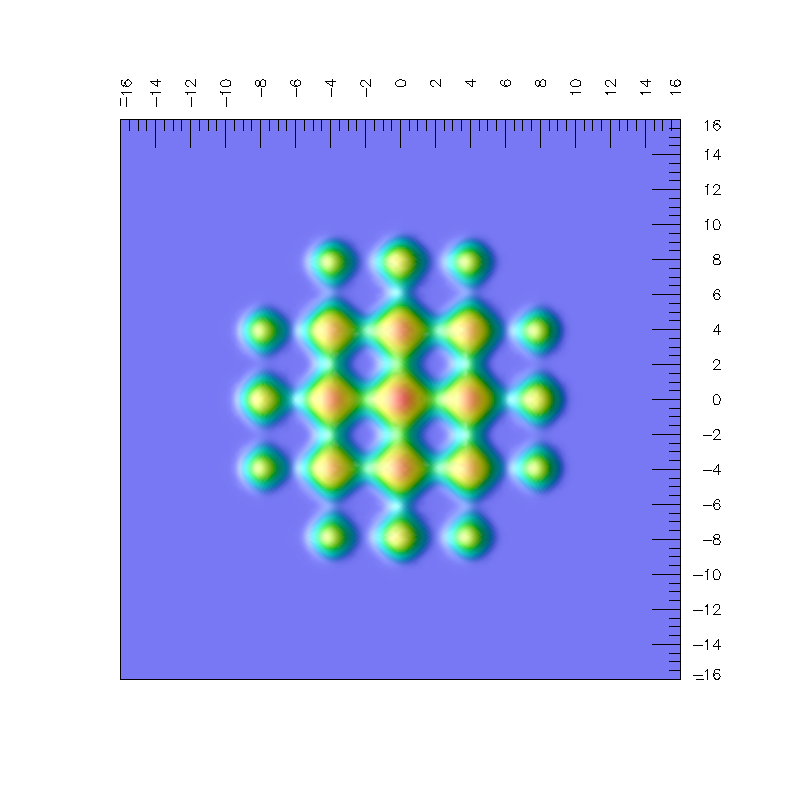}
        \includegraphics[width=0.3\linewidth]{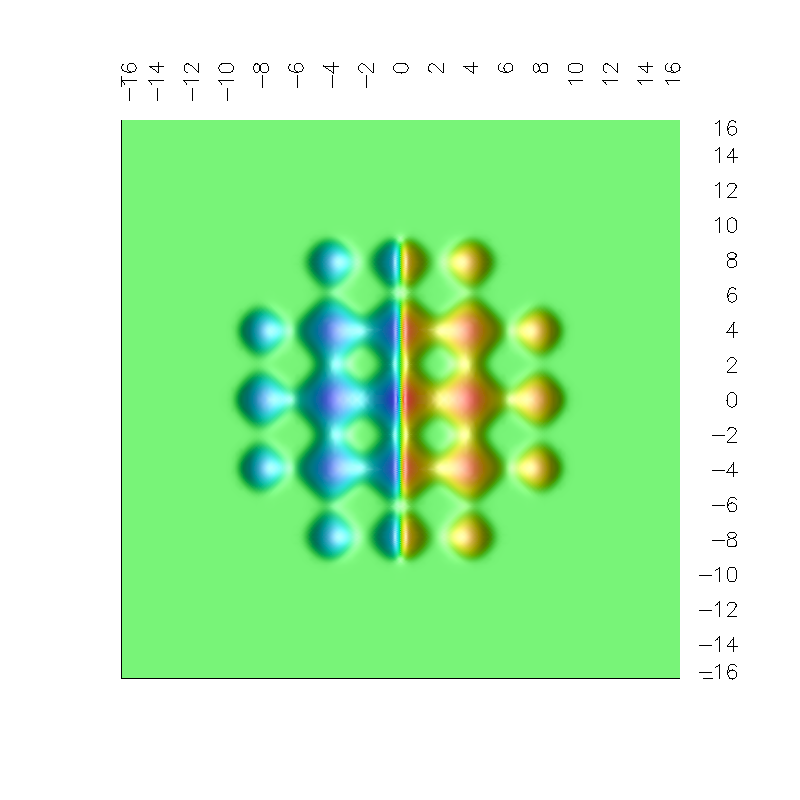}
        \includegraphics[width=0.3\linewidth]{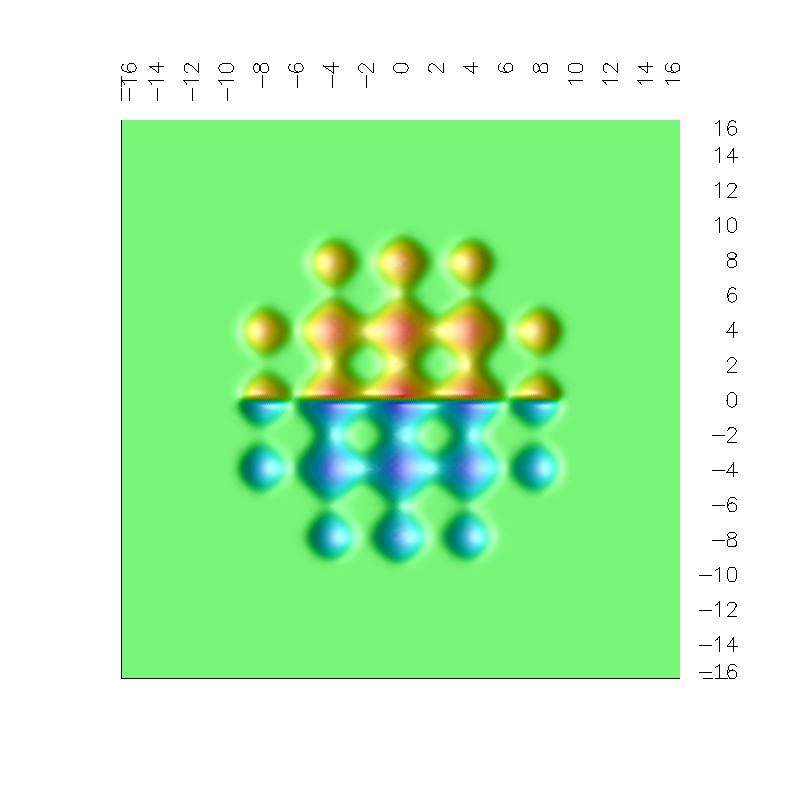}
        \caption{Simultaneously solving the desired states of multi-mesh in Example 2. Top: meshes of ground state, $x$-direction excited state and $y$-direction excited state; Bottom: solutions of ground state, $x$-direction excited state and $y$-direction excited state.}
        \label{mul_excited_meshes}
    \end{figure}

    \begin{table}[htbp]
        \centering 
        \vspace{0.5em}
        \caption{Comparison on desired states of Example 2 using the single-mesh and multi-mesh under the same adaptive tolerance $10^{-4}$.}
        \vspace{0.5em}
        \label{tab:comparison2}
        \begin{tabular}{lllllll} 
            \toprule
            {} & $E_{00}$ & $E_{10}$ & $E_{01}$ & $N_{sol1}$ & $N_{sol2}$ & $N_{sol3}$ \\
            \midrule
            single-mesh & 55.443 & 57.146 & 57.146 & 98,889 & 98,889 & 98,889 \\
            multi-mesh & 55.443 & 57.146 & 57.146 & 35,961 & 61,249 & 61,249 \\
            \bottomrule
        \end{tabular}
    \end{table}

    The meshes of the single-mesh method and the multi-mesh method are shown in Figure \ref{single_excited_mesh} and \ref{mul_excited_meshes}, respectively.  The top of Figure \ref{mul_excited_meshes} from left to right is the meshes of ground state, $x$-direction excited state and $y$-direction excited state, and the bottom is the solution of ground state, $x$-direction excited state and $y$-direction excited state, respectively. From these meshes, we can see that the multi-mesh method can adaptively refine the mesh according to the characteristics of different states.

    In Table \ref{tab:comparison2}, $E_{00}$, $E_{10}$, and $E_{01}$ denote the energies of the ground state, the excited state in the $x$-direction, and the excited state in the $y$-direction, respectively. The quantities $N_{\text{sol}1}$, $N_{\text{sol}2}$, and $N_{\text{sol}3}$ indicate the numbers of degrees of freedom required to compute these three states. As shown in the table, to attain the same level of numerical accuracy as the single-mesh method, which uses 98,889 degrees of freedom, the multi-mesh method requires only 35,961 degrees of freedom for the ground state and 61,249 degrees of freedom for the excited states. It is worth noting that although the two meshes used for the excited states are not identical, they happen to have the same number of degrees of freedom.  These results demonstrate that the multi-mesh method is more efficient than the single-mesh method for achieving comparable numerical accuracy in computing both ground and excited states.


    \subsection*{Example 3}
    In the last example, we take the two-dimensional condensate \cite{Caliari2009,Kuang2018} into account. We use the multi-mesh adaptive finite element method proposed to calculate the ground state of the coupled Gross-Pitaevskii system and take the potential functions scaled and off centered harmonic potentials:
    $$U_i = \frac{1}{2} \left[ (\omega_{1i}(x - x_i))^2 + (\omega_{2i}(y - y_i))^2 \right], i = 1,2,$$
    with
    $$\omega_{11} = \omega_{21} = \pi, \, x_1 = y_1 = 0, \, \omega_{12} = \omega_{22} = 3\pi, \, x_2 = 0.19, \, y_2 = 0,$$
    and the intra-species coupling constants
    $$v_{11} = 1.3 \times 10^{-6}, \quad v_{12} = v_{21} = 1.0 \times 10^{-6}, \quad v_{22} = 1.3 \times 10^{-11}.$$

\begin{figure}[htbp]
        \centering
        \begin{minipage}{0.23\linewidth}
            \centering
            \includegraphics[width=\linewidth]{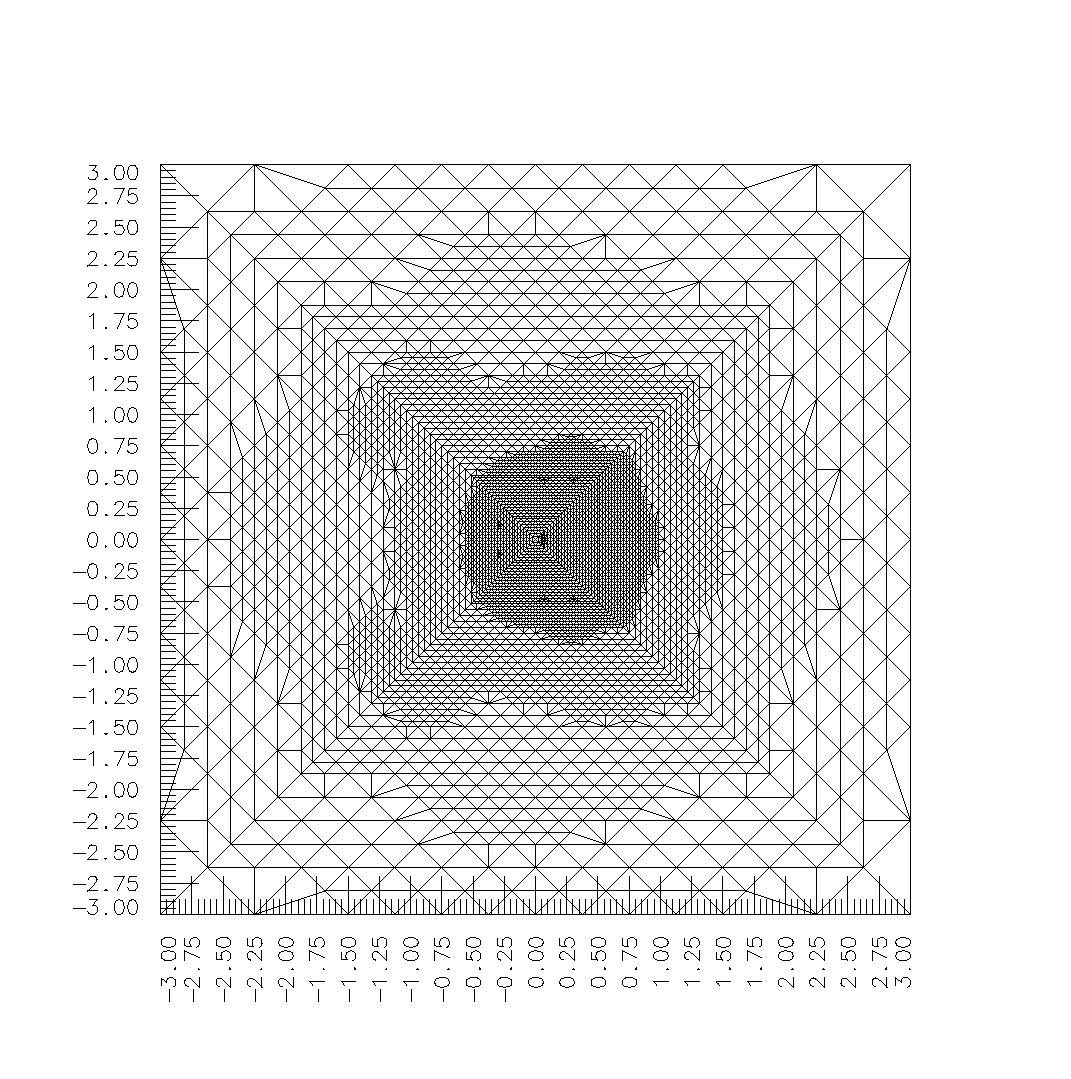}
        \end{minipage}
        \begin{minipage}{0.23\linewidth}
            \centering
            \makebox[0.2\linewidth][c]{\vphantom{\includegraphics[width=0.2\linewidth]{Example3_Figure/single_mesh.png}}}
        \end{minipage}
        \begin{minipage}{0.23\linewidth}
            \centering
            \includegraphics[width=\linewidth]{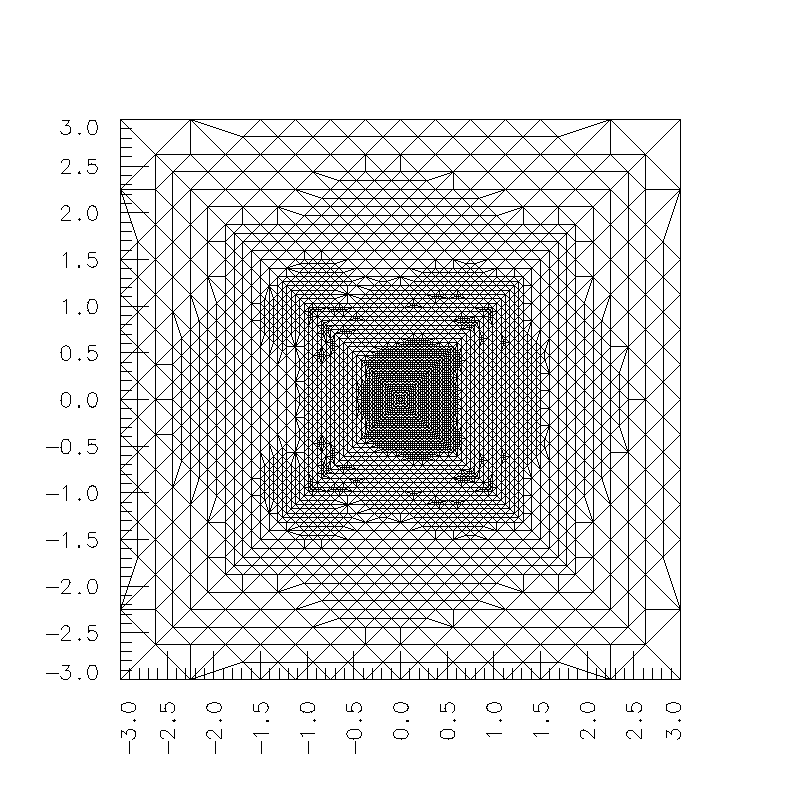}
        \end{minipage}
        \begin{minipage}{0.23\linewidth}
            \centering
            \includegraphics[width=\linewidth]{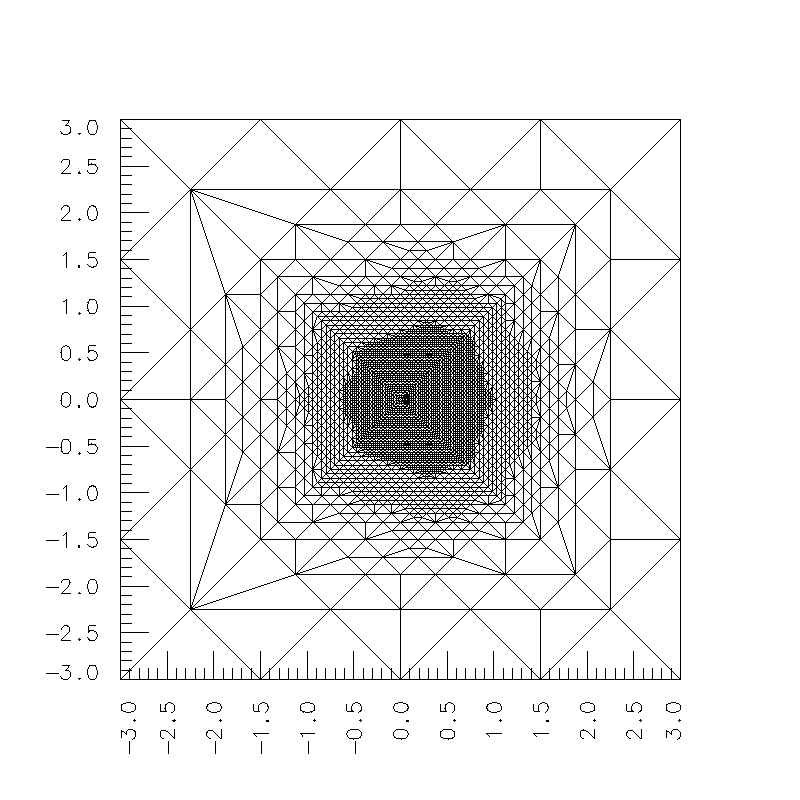}
        \end{minipage}

        \vspace{1ex}

        \begin{minipage}{0.23\linewidth}
            \centering
            \includegraphics[width=\linewidth]{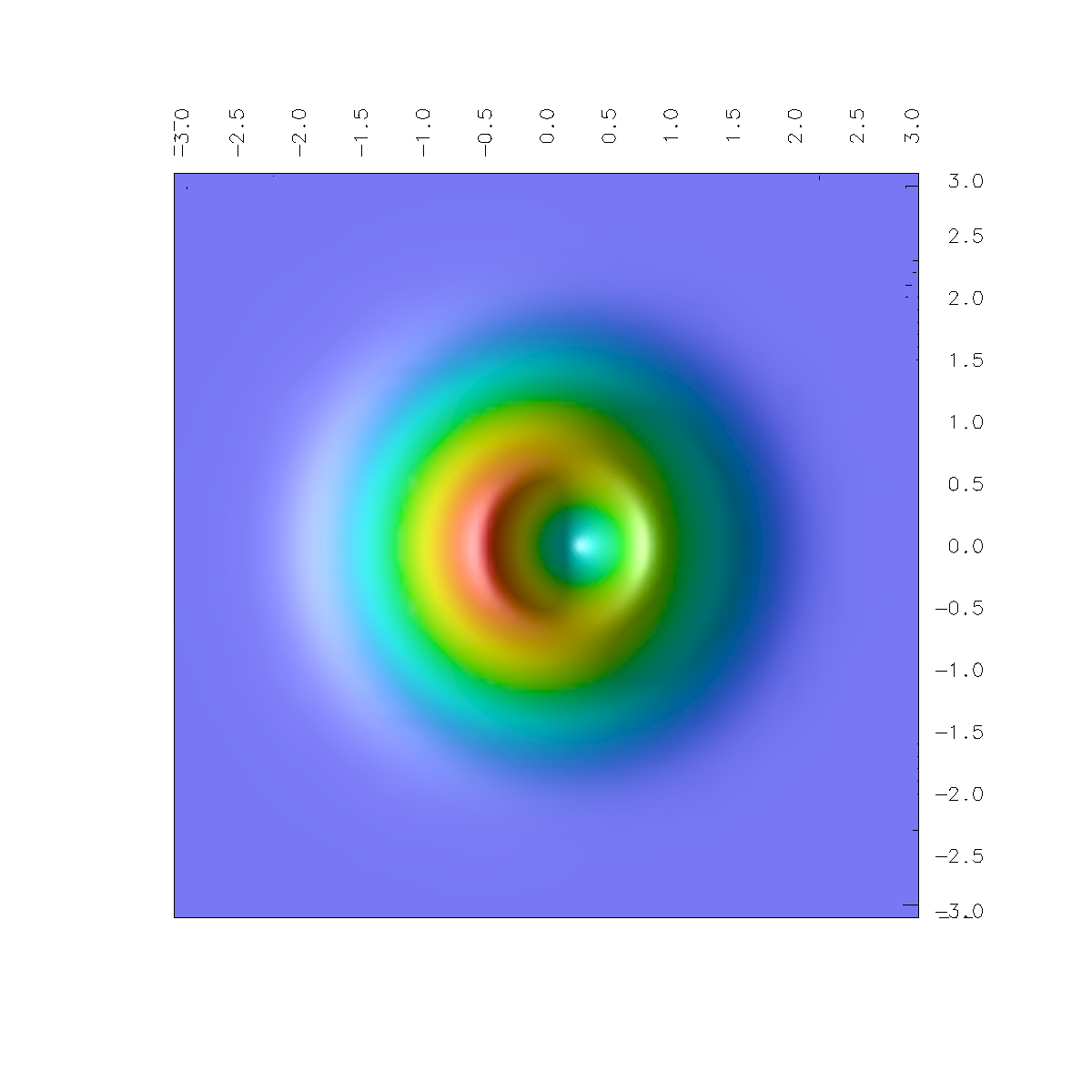}
        \end{minipage}
        \begin{minipage}{0.23\linewidth}
            \centering
            \includegraphics[width=\linewidth]{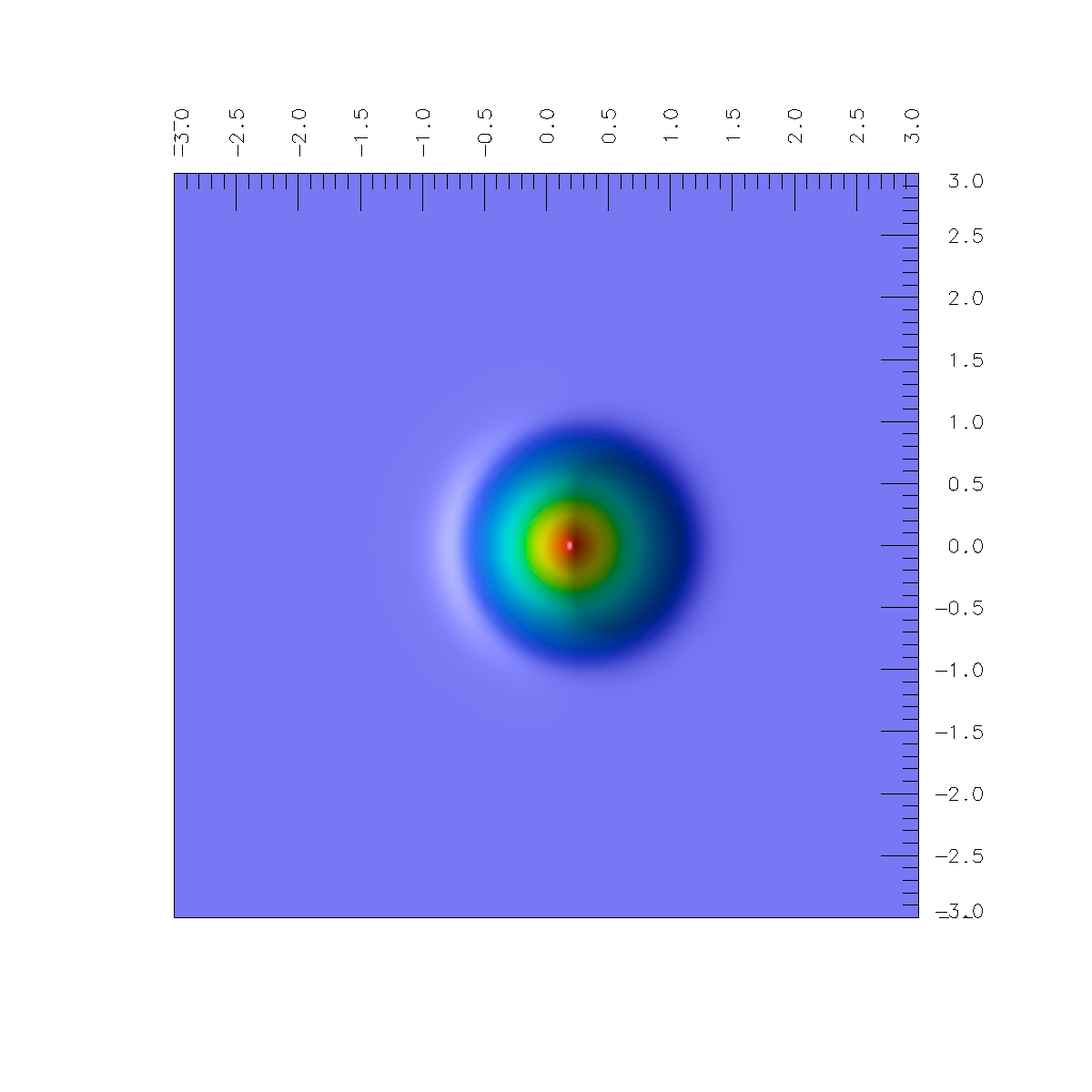}
        \end{minipage}
        \begin{minipage}{0.23\linewidth}
            \centering
            \includegraphics[width=\linewidth]{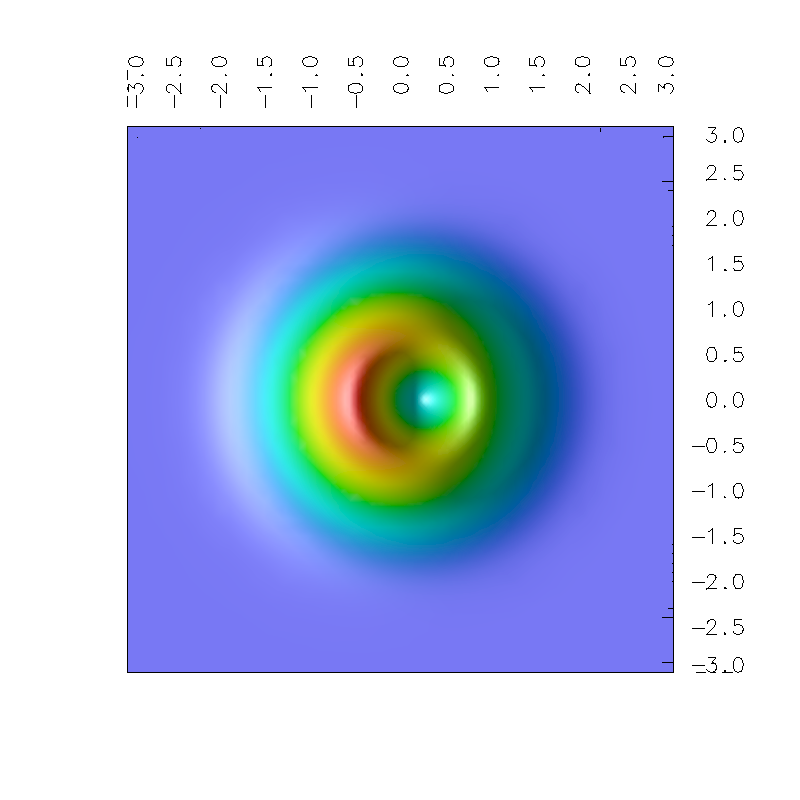}
        \end{minipage}
        \begin{minipage}{0.23\linewidth}
            \centering
            \includegraphics[width=\linewidth]{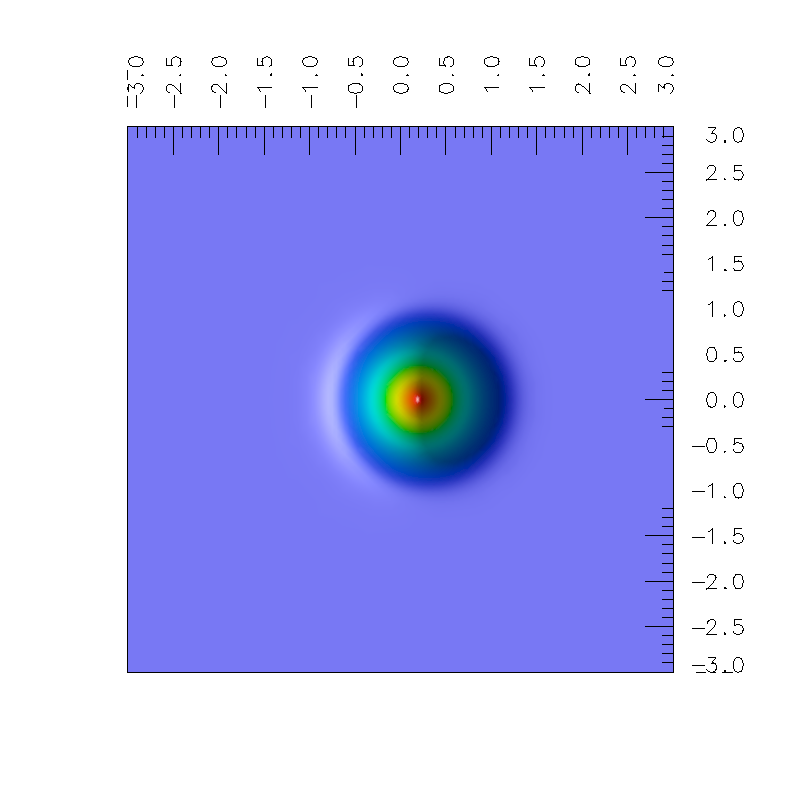}
        \end{minipage}

        \caption{Meshes and solutions of single-mesh method and multi-mesh method in ground state of Example 3. Top left: single-mesh, top right: multi-mesh.}
        \label{mulgpe}
    \end{figure}

    We set the same number of the particles $N_1=N_2=10^7.$ We conducted the numerical experiment on region $[-3,3]\times[-3,3]$ and compared the results obtained from the multi-mesh adaptive finite element method with those obtained from the single-mesh method. The experimental results are shown in Figure \ref{mulgpe}, the three images above are: on the left are the mesh for solving the coupled Gross-Pitaevskii system using the single-mesh adaptive finite element method, and the two on the right correspond to the meshes for solving two solutions in the coupled Gross-Pitaevskii system using the multi-mesh adaptive finite element method. The four images below are: two solutions in coupled Gross-Pitaevskii system, the left two solutions using a single-mesh method, and the right two solutions using the multi-mesh method.

    In Table \ref{tab:comparison3}, $N_{sol1}$ and $N_{sol2}$ represent the number of degrees of freedom used for the first solution and second solution, respectively. We used the solution under a uniform mesh with degrees of freedom of 525,313 as a reference value, and the total energy and total chemical potential were calculated as 169,685,071.853351 and 192,904,970.150445, respectively.  When using the single grid adaptive finite element method, we need to use 15793 degrees of freedom to achieve numerical accuracy, while using the multi-mesh method reducing the degrees of freedom by 11 percent and 29 percent, respectively. The experimental results indicate that our proposed method can achieve the desired numerical accuracy with less computational cost.

    \begin{table}[htbp]
        \centering 
        \vspace{0.5em}
        \caption{Comparison of numerical results obtained with single-mesh and multi-mesh methods for the multi-component GPEs in Example 3. The referenced total energy is $E_{ref} = 169,685,071.853351$ and the referenced total chemical potential is $\mu_{ref} = 192,904,970.150445$.}
        \vspace{0.5em}
        \label{tab:comparison3}
        \begin{tabular}{lllll} 
            \toprule
            {} & $E_{total}$ & $\mu_{total}$ & $N_{sol1}$ & $N_{sol2}$ \\
            \midrule
            single-mesh & 169,685,247.085354 & 192,905,174.229224 & 15,793 & 15,793 \\
            multi-mesh & 169,685,269.297038 & 192,905,230.772153 & 14,001 & 11,161 \\
            \bottomrule
        \end{tabular}
    \end{table}

    \section{Conclusion}
    
    In this paper, we introduce a multi-mesh adaptive finite element method for solving the GPE, with the goal of improving computational efficiency. Compared with uniform mesh refinement, adaptive mesh refinement achieves the desired accuracy using significantly fewer elements by concentrating resolution only where needed. Moreover, to accurately compute multiple states such as both ground and excited states, we extend this approach to a multi-mesh framework, in which each state is solved on its own adaptively refined mesh, thereby further reducing the computational cost. The proposed strategy also applies naturally to multi-component GPEs, where different components may exhibit distinct spatial variations. Numerical experiments demonstrate the effectiveness and accuracy of the proposed multi-mesh adaptive finite element method.

    \section*{Acknowledgement}
    The work of Y. Kuang is supported by the National Natural Science Foundation of China (Nos. 12201130 and 12326362), and fund from the Key Laboratory of Mathematical Modelling and High Performance Computing of Air Vehicles (NUAA) (No. 202302). The work of Z. Hu is supported by the National Natural Science Foundation of China (No. 12171240), and fund from the Key Laboratory of Mathematical Modelling and High Performance Computing of Air Vehicles (NUAA) (No. 202302).
    \bibliographystyle{plain}
    \bibliography{references}
    \end{document}